\documentclass[twoside,onecolumn]{elsart1p}
\usepackage{amsmath,amssymb,latexsym,theorem,bbm,epic,eepic}

\newcommand{\SC}{\scriptstyle}
\newcommand{\bone}{\mathbbm{1}}
\newcommand{\stoch}{\stackrel{\SC\mathsf{P}}{\longrightarrow}}
\newcommand{\distr}{\stackrel{\SC{\mathcal D}}{\longrightarrow}}
\newcommand{\qmean}{\stackrel{{\SC{\mathsf L}_2}}{\longrightarrow}}

\def\cov{\mathop{\hbox {\sf Cov}}}
\def\var{\mathop{\hbox {\sf Var}}}
\def\sign{\mathop{\hbox {\rm sign}}}
\def\diag{\mathop{\hbox {\rm diag}}}
 
\numberwithin{equation}{section}

\newcommand{\proofend}{\hfill$\square$}

\theorembodyfont{\em}
\newtheorem{Lem}{Lemma}[section]
\newtheorem{Thm}[Lem]{Theorem}
\newtheorem{Pro}[Lem]{Proposition}

\newtheorem{Cor}[Lem]{Corollary}
\theorembodyfont{\rm}
\newtheorem{Rem}[Lem]{Remark}

\begin{document}
\begin{frontmatter}

\title{On the least squares estimator in  a
  nearly unstable sequence of stationary spatial AR models}

\author[gp]{S\'andor Baran\corauthref{bscorr}},
\ead{barans@inf.unideb.hu}
\author[gp]{Gyula Pap}
\address[gp]{Faculty of Informatics, University of
             Debrecen, Hungary}
\corauth[bscorr]{Corresponding author. Address:
                 Faculty of Informatics, University of
                 Debrecen, H-4010 Debrecen, P.O. Box 12, Hungary. Tel:
                 +36 52 512900, fax: +36 52 416857.}

\begin{abstract}
A nearly unstable sequence of stationary spatial autoregressive processes is
 investigated, when the sum of the absolute values of the
 autoregressive coefficients tends to one. 
It is shown that after an appropriate norming the least squares estimator for 
 these coefficients has a normal limit distribution. If none of the
 parameters equals zero than the typical rate of convergence is $n$.

\medskip
\noindent
\textit{AMS} 2000 \textit{subject classifications.} Primary 62M10;
Secondary 62F12
\end{abstract}

\end{frontmatter}

\section{Introduction}
   \label{sec:sec1}
Spatial autoregressive models have a great importance 
in many different fields of science such as geography, geology,
biology and agriculture, see e.g. \cite{br2} for a detailed
discussion, where the authors considered a general unilateral model
having the form  
 \begin{equation}
    \label{eq:eq1.1}
  X_{k_,\ell}
  =\sum_{i=0}^{p_1}\sum_{j=0}^{p_2}\alpha_{i,j}X_{k-i,\ell-j}
   +\varepsilon_{k,\ell}, 
  \qquad\alpha_{0,0}=0.
\end{equation}
A particular case of the model \eqref{eq:eq1.1} is the
 so-called doubly geometric spatial autoregressive model
\begin{equation*}
  X_{k,\ell}
  =\alpha X_{k-1,\ell}+\beta X_{k,\ell-1}-\alpha\beta X_{k-1,\ell-1}
   +\varepsilon_{k,\ell},
 \end{equation*}
 introduced by Martin \cite{martin1}. In fact, this is the simplest spatial
 model, since its nice product structure ensures that it can be
 considered as some kind of combination of two 
 autoregressive processes on the line, and several properties can be derived
 by the analogy of one-dimensional autoregressive processes. The
 doubly geometric model was the first one for which the nearly
 unstability has been 
studied. Bhattacharyya {\em et al.\/} \cite{bhat2} showed that in the
case when a sequence of stable models with \
 $\alpha_n\to1$, $\beta_n\to1$ \ was considered, in contrast to the AR(1) model,
 the sequence of Gauss-Newton estimators \
 $(\widehat\alpha_n,\widehat\beta_n)$ \ of 
 \ $(\alpha_n,\beta_n)$ \ were asymptotically normal, namely, \
\begin{equation*}
n^{3/2}\begin{pmatrix}
       \widehat\alpha_n-\alpha_n \\ \widehat\beta_n-\beta_n
       \end{pmatrix}
  \distr{\mathcal N}(0,\Sigma)
\end{equation*}
 with some covariance matrix \ $\Sigma$. \

The doubly geometric model has several
applications. Jain \cite{jain} used it in the study of image processing, 
Martin \cite{martin2}, Cullis and Gleeson \cite{cg}, Basu and Reinsel
\cite{br1} in agricultural trials, while Tj{\o}stheim
\cite{tj1} in digital filtering. 

In the present paper we study another special case of the model
\eqref{eq:eq1.1}. We consider the spatial autoregressive process \
$\{X_{k,\ell}:k,\ell\in{\mathbb Z}\}$ \ which is a
 solution of the spatial stochastic difference equation
 \begin{equation}\label{eq:eq1.2}
  X_{k,\ell}=\alpha X_{k-1,\ell}+\beta X_{k,\ell-1}+\varepsilon_{k,\ell}
 \end{equation}
with parameters \ $(\alpha,\beta)\in {\mathbb R}^2$. \ 
This model is stable (i.e. has a stationary solution) in  case \
$|\alpha|+|\beta |<1$ \
(see \cite{br2}), and unstable if \ $|\alpha|+|\beta |=1$. \
In a recent paper Paulauskas \cite{paul} determined the exact
asymptotic behavior of the variances of a nonstationary solution of
 \eqref{eq:eq1.2}  with \ $X_{k,\ell}=0$ \ for \
 $k+\ell\leq0$, \ while Baran {\em et
  al.\/} \cite{bpz3} in the same model
clarified the asymptotic properties of the least
 squares estimator (LSE)
of \ $(\alpha,\beta )$ \  both in stable and unstable cases.

We remark, that in  case \ $|\alpha|+|\beta |<1$, if \
$\{\varepsilon_{k,\ell}:k,\ell\in{\mathbb Z}\}$ \ are independent and
identically distributed random  variables,  a stationary solution
can be given by
 \begin{equation}\label{eq:eq1.3}
  X_{k,\ell }
  =\sum_{(i,j)\in\, U_{k,\ell}}\binom{k+\ell -i-j}{k-i}
    \alpha^{k-i}\beta^{\ell -j}\varepsilon_{i,j},
 \end{equation}
 where \ $U_{k,\ell}:=\{(i,j)\in{\mathbb Z}^2:\text{$i\leq k$ and $j\leq
   \ell$}\}$ \ and
 the convergence of the series is understood in ${\mathsf L}_2$-sense.

We are interested in the asymptotic behaviour of the stationary
solution of \eqref{eq:eq1.2} in the case when the parameters approach
the boundary \ $|\alpha |+|\beta |=1$. \ In order to determine the
appropriate speed of parameters one may use the idea of Chan and
Wei \cite{cw} and consider the order of
\begin{equation*}
{\mathbb I}_n:={\mathsf E} \left(\sum_{(k,\ell)\in H_n}
               \begin{pmatrix}
                 \big(X_{k-1,\ell}\big )^2 & X_{k-1,\ell}
                        X_{k,\ell-1} \\
                 X_{k-1,\ell}X_{k,\ell-1} & 
                 \big(X_{k,\ell-1} \big)^2
               \end{pmatrix} \right)
\end{equation*}
that is exactly the observed Fisher information matrix about \
$(\alpha,\beta)$ \ when the innovations \ $\varepsilon _{k,\ell}$ \
are normally distributed and the process is observed on a set \
$H_n\subset {\mathbb Z}^2, \ n\in{\mathbb N}$. \ From Theorem 1.1 of
\cite{bpz3} we obtain that
\begin{equation*}
{\mathbb I}_n \sim \begin{cases}
                      n^2\sigma^2_{\alpha,\beta}\Gamma_{\alpha,\beta},
                      & \text{if \ $|\alpha|+|\beta 
                        |<1$}, \\[2 mm]
                      n^{5/2}\sigma^2_{\alpha }\Psi_{\alpha,\beta}, &
                      \text{if \ $|\alpha|+|\beta|=1, \ 0<|\alpha|
                        <1$}, \\[2 mm] 
                       n^3(4/3){\mathcal I}, &
                      \text{if \ $|\alpha|+|\beta|=1, \ |\alpha|\in
                        \{0,1\}$},  
                    \end{cases}
\end{equation*}
where 
\begin{equation*}
\Gamma_{\alpha,\beta}
    :=2
      \begin{pmatrix}
         1 & -\varrho_{\alpha,\beta}\\
        -\varrho_{\alpha,\beta} & 1
      \end{pmatrix}, \qquad \qquad
\Psi _{\alpha,\beta}:=\begin{pmatrix}
        1 & \sign (\alpha\beta)\\
        \sign(\alpha\beta) & 1
      \end{pmatrix},
\end{equation*}
${\mathcal I}$ \ denotes the two-by-two unit matrix and
 \begin{align*}
  \sigma_{\alpha,\beta}^2
  &:=\big((1+\alpha+\beta)(1+\alpha-\beta)(1-\alpha+\beta)
                               (1-\alpha-\beta)\big)^{-1/2},\\[2mm]
   \varrho_{\alpha,\beta} 
   &:=\begin{cases}
          \displaystyle\frac{(1-\alpha^2-\beta^2)\sigma_{\alpha,\beta}^2-1}
                         {2\alpha\beta\sigma_{\alpha,\beta}^2},
         &\text{if \ $\alpha\beta\not=0$,}\\
         0 &\text{otherwise,}  
       \end{cases} \\[2mm]
 \sigma_\alpha^2&:=\frac{2^{9/2}}{15\sqrt{\pi|\alpha|(1-|\alpha|)}}.
\end{align*}
Now, let \ $\alpha _n:=\alpha -\gamma/a_n, \ \beta_n:=\beta
-\delta/a_n, \ |\alpha|+|\beta|=1, \ |\alpha_n|+|\beta_n|<1$. \ As
nonstationary behaviour of \ 
$X_{k,\ell}$ \ becomes dominant when \ $(\alpha_n, \beta_n)$ \ is near
the border, a reasonable choice for the sequence \ $a_n$ \ should
retain the order of \ ${\mathbb I}_n$ \ to be \ $n^{5/2}$ \ if \
$0<|\alpha|<1$ \ and \ $n^3$ \ if \ $|\alpha|\in \{0,1\}$. \
Since we have \ $\sigma _{\alpha_n,\beta_n}^2\sim a_n^{1/2}$ \  for
$0<|\alpha|<1$ \ and  \ $\sigma _{\alpha_n,\beta_n}^2\sim a_n$ \
for \ $|\alpha|\in \{0,1\}$ \ while  \  $\varrho
_{\alpha_n,\beta_n}\sim const$ \ in both cases, the above
consideration yields \ $a_n=n$. 

In what follows we consider a nearly unstable sequence of stationary processes,
 i.e. for each \ $n\in{\mathbb N}$, \ we take a stationary solution
 \ $\{X^{(n)}_{k,\ell}:k,\ell\in{\mathbb Z}\}$ \ of equation
 \eqref{eq:eq1.2} with
 parameters \ $(\alpha_n, \beta_n)$ \ defined as
 \begin{equation}\label{eq:eq1.4}
  \alpha _n:=\alpha -\frac {\gamma _n}n, \qquad \beta _n:=\beta -
\frac {\delta _n}n, \quad \qquad |\alpha _n|+|\beta_n|<1,
 \end{equation}
where \  $0\leq |\alpha| \leq1, \ |\beta|=1-|\alpha |$ \ and \
$\gamma_n\to\gamma,  \ \delta_n\to\delta$ \ as \ $n\to\infty$, \
$(\gamma,\delta)\in {\mathbb R}^2$. \ We remark that in an earlier
paper \cite{bpz2} the authors considered a similar sequence of
stationary processes  where the autoregressive parameters were equal and
their sum converged to \ $1$. 

For a set \ $H\subset{\mathbb Z}^2$, \ the LSE  \
$(\widehat\alpha^{(n)}_H,\widehat\beta^{(n)}_H)$ \ of \
$(\alpha_n,\beta_n)$ \ based on the observations  \
$\{X^{(n)}_{k,\ell}:(k,\ell)\in H\}$ \ has the form
\begin{equation*}
\begin{pmatrix}\widehat\alpha_H^{(n)}\\
                   \widehat\beta_H^{(n)}\end{pmatrix} =
\left(\sum_{(k,\ell)\in H}
               \begin{pmatrix}
                 \big (X_{k-1,\ell}^{(n)}\big )^2 & X_{k-1,\ell}^{(n)}
                        X_{k,\ell-1}^{(n)} \\
                 X_{k-1,\ell}^{(n)}X_{k,\ell-1}^{(n)} & \big
                 (X_{k,\ell-1}^{(n)}\big )^2
               \end{pmatrix}\right )^{-1}
\sum_{(k,\ell)\in H}
       \begin{pmatrix}
        X_{k-1,\ell}^{(n)}X_{k,\ell}^{(n)} \\
        X_{k,\ell-1}^{(n)}X_{k,\ell}^{(n)}
       \end{pmatrix}.
 \end{equation*}
Consider the triangles \
 $T_{k,\ell}
  :=\{(i,j)\in{\mathbb Z}^2:i+j\geq1,\,i\leq
  k\,\,\text{and}\,\,j\leq\ell\}$ \
 for \ $k,\ell\in{\mathbb Z}$. \
Note that \ $T_{k,\ell}=\emptyset$ \ if \ $k+\ell\leq0$.

\begin{Thm}
   \label{main} \
For each \ $n\in{\mathbb N}$, \ let \ $\{X^{(n)}_{k,\ell}:k,\ell\in{\mathbb
  N}\}$ \ be a stationary  solution of equation \eqref{eq:eq1.2} with
parameters \ $(\alpha_n, \beta_n)$ \ given by
 \eqref{eq:eq1.4}, and with independent and identically distributed random
 variables \ $\{\varepsilon^{(n)}_{k,\ell}:k,\ell\in{\mathbb Z}\}$ \ such that
 \ ${\mathbb E}\varepsilon^{(n)}_{0,0}=0$, \
 $\var\varepsilon^{(n)}_{0,0}=1$ \ and \
 $M:=\sup_{n\in{\mathbb N}}{\mathbb
   E}\big|\varepsilon^{(n)}_{0,0}\big|^8<\infty$. \
Let \ $(k_n)$ \ and \ $(\ell_n)$ \ be sequences of integers such that \
 $k_n+\ell_n\to\infty$ \ as \ $n\to\infty$. \

If \ $0<|\alpha |<1, \ |\beta|=1-|\alpha |$ \ and
\begin{equation}
   \label{eq:eq1.5}
\lim _{n\to\infty }(k_n+\ell _n)n^{-1/2}\big (|\gamma _n|+|\delta
_n|\big )^{1/2}=\infty
\end{equation}
holds then
\begin{equation*}
(k_n+\ell _n)
  \begin{pmatrix}\widehat\alpha_{T_{k_n,\ell _n}}-\alpha_n\\
                   \widehat\beta_{T_{k_n,\ell
                       _n}}-\beta_n\end{pmatrix}
\distr {\mathcal N}_2\Big (0,|\alpha||\beta|\bar \Psi_{\alpha,\beta}\Big )
\end{equation*}
as \ $n\to\infty$, \ where \ $\bar \Psi_{\alpha,\beta}$ \ denotes the
adjoint matrix of \ $ \Psi _{\alpha, \beta}$.

If \ $|\alpha |\in \{0,1\},  \ |\beta|=1-|\alpha |$ \ and
\begin{equation}
   \label{eq:eq1.6}
\lim _{n\to\infty }(k_n+\ell _n)n^{-1}\big |\gamma _n^2-\delta
_n^2\big |^{1/2}=\infty
\end{equation}
holds then let
\begin{equation*}
[-\infty,\infty]\ni \omega:=\lim_{n\to\infty}\omega _n, \qquad \qquad
\omega_n:=\alpha 
\frac{\gamma_n}{\delta_n}+\beta  \frac{\delta_n}{\gamma_n}.
\end{equation*}
If \ $|\omega | >1$ \ then
\begin{equation*}
 (k_n+\ell _n)n^{1/2}\big |\gamma _n^2-\delta _n^2\big |^{-1/4}
  \begin{pmatrix}\widehat\alpha_{T_{k_n,\ell _n}}-\alpha_n\\
                   \widehat\beta_{T_{k_n,\ell _n}}-\beta_n\end{pmatrix}
  \distr {\mathcal N}_2 \Big (0,\Theta _{\alpha,\beta,\omega}^{-1}\Big)
\end{equation*}
as \ $n\to\infty$, \ where
\begin{equation*}
 \Theta _{\alpha,\beta,\omega}:=\frac 14 \begin{pmatrix} 1& \theta
   (\alpha,\beta, \omega) \\
   \theta (\alpha,\beta, \omega)&
                           1\end{pmatrix}
\end{equation*}
with
\begin{equation*}
\theta (\alpha,\beta, \omega):= \begin{cases}
                                \frac {-(\alpha +\beta )\sign(\omega)}{|\omega|
                                  +\sqrt {\omega ^2-1}}& \text{if \
                                  $|\omega |<\infty $,} \\
                                 0& \text {if \ $|\omega| =\infty$.}
                                \end{cases}
\end{equation*}
\end{Thm}

\begin{Rem}
  \label{rem:rem1} 
Obviously, \ $|\omega _n|>1$, \ so \ $|\omega|\geq1$.\ Condition \
$|\omega |>1$ \ in Theorem \ref{main} is needed to ensure the regularity of \  
$\Theta _{\alpha,\beta,\omega}$. \ However, this condition can be omitted and
using similar arguments as in the proof of the second statement of
Theorem \ref{main}, one can easily show that if \ $|\alpha |\in
\{0,1\},  \ |\beta|=1-|\alpha |$ \ and \eqref{eq:eq1.6} holds then
\begin{equation*}
 (k_n+\ell _n)n^{1/2}\big |\gamma _n^2-\delta _n^2\big
 |^{-1/4}\Theta_{\alpha,\beta,\omega_n}^{1/2} 
  \begin{pmatrix}\widehat\alpha_{T_{k_n,\ell _n}}-\alpha_n\\
                   \widehat\beta_{T_{k_n,\ell _n}}-\beta_n\end{pmatrix}
  \distr {\mathcal N}_2 \Big (0,{\mathcal I}\Big),
\end{equation*}
where 
\ $\Theta_{\alpha,\beta,\omega_n}^{1/2}$ \ denotes the symmetric  positive
semidefinite square root of \ $\Theta_{\alpha,\beta,\omega_n}$. 
\end{Rem}

\begin{Rem}
  \label{rem:rem2}
Theorem \ref{main} shows that in the typical case \ $k_n=\ell_n=n$ \ and \
$\gamma_n=\gamma\ne 0$, \ $\delta_n=\delta\ne 0$ \ if \ $0<|\alpha
|<\infty, \ |\beta|=1-|\alpha |$ \ then the rate of
convergence is \ $n$.
\end{Rem}

We may suppose that \ $(k_n+\ell_n)$ \ is monotone increasing.
Observe, that \ $\big (\widehat\alpha^{(n)}_{T_{k_n,\ell_n}},
\widehat\beta^{(n)}_{T_{k_n,\ell_n}}\big )$ \ and \
 $\big (\widehat\alpha^{(n)}_{T_{\widetilde
     k_n,\widetilde\ell_n}},\widehat\beta^{(n)}_{T_{\widetilde
     k_n,\widetilde\ell_n}}\big )$ \ have the same distribution, where \
 $\widetilde k_n:=[(k_n+\ell_n)/2]$ \ and \
 $\widetilde\ell_n:=[(k_n+\ell_n+1)/2]$. \
As \ $\widetilde k_n+\widetilde \ell_n=k_n+\ell_n$, \ in Theorem \ref{main}
we may substitute \ $(\widetilde k_n,\widetilde\ell_n)$ \ for \
$(k_n,\ell_n)$. \
The sequence \ $(\widetilde k_n,\widetilde\ell_n)$ \ can be embedded into
the sequence \ $(k'_n,\ell'_n)$, \ where \ $k'_n:=[n/2]$ \ and \
$\ell'_n:=[(n+1)/2]$, \ namely, \ $k'_{q_n}=\widetilde k_n$ \ and \
$\ell'_{q_n}=\widetilde\ell_n$ \ with \ $q_n:=\widetilde
k_n+\widetilde\ell_n$.  \
Clearly \ $k'_n+\ell'_n=n$. \ Consider the sequence \ $(r_n)$ \ defined by
\ $r_n:=k$ \ for \ $q_k\leq n<q_{k+1}$. \ Then \ $r_{q_n}=n$, \ and conditions
\eqref{eq:eq1.5} and \eqref{eq:eq1.6} can be replaced by
\begin{equation}
   \label{eq:eq1.7}
\lim _{n\to\infty }nr_n^{-1/2}\big (|\gamma _{r_n}|+|\delta
_{r_n}|\big )^{1/2}=\infty
\end{equation}
and
\begin{equation}
   \label{eq:eq1.8}
\lim _{n\to\infty }nr_n^{-1}\big|\gamma _{r_n}^2-\delta _{r_n}^2\big
|^{1/2}=\infty,
\end{equation}
respectively.

Thus, to prove Theorem \ref{main} it suffices to show that if
\ $0<|\alpha |<1,  \ |\beta|=1-|\alpha |$ \ and \eqref{eq:eq1.7} holds then
\begin{equation*}
n \begin{pmatrix}\widehat\alpha_{T_{[n/2],[(n+1)/2]}}-\alpha_{r_n}\\
                   \widehat\beta_{T_{[n/2],[(n+1)/2]}}-\beta_{r_n}\end{pmatrix}
  \distr{\mathcal N}_2\Big (0,|\alpha||\beta|\bar \Psi_{\alpha,\beta}\Big ),
\end{equation*}
while in the case  \ $|\alpha |\in \{0,1\}, \ |\beta|=1-|\alpha |$ ,
\ $|\omega | >1$ \ and \eqref{eq:eq1.8} holds
we have
\begin{equation*}
 nr_n^{1/2}\big |\gamma _{r_n}^2-\delta _{r_n}^2\big |^{-1/4}
  \begin{pmatrix}\widehat\alpha_{T_{[n/2],[(n+1)/2]}}-\alpha_{r_n}\\
                   \widehat\beta_{T_{[n/2],[(n+1)/2]}}-\beta_{r_n}\end{pmatrix}
  \distr {\mathcal N}_2 \Big (0,\Theta
  _{\alpha,\beta,\omega}^{-1}\Big).
\end{equation*}

We remark that conditions \eqref{eq:eq1.5} and \eqref{eq:eq1.7} are
exactly the same as conditions (4) and (5) of \cite{bpz2}, respectively.

To simplify notation we assume \ $k_n=[n/2]$, \ $\ell _n=[(n+1)/2]$
\ and \ $(r_n)$ \ is a monotone increasing sequence of positive
integers. One can write 
\begin{equation*}
\begin{pmatrix}\widehat\alpha_{T_{k_n,\ell _n}}-\alpha_{r_n}\\
                   \widehat\beta_{T_{k_n,\ell _n}}-\beta_{r_n}\end{pmatrix}
    =B_n^{-1}A_n,
\end{equation*}
with
 \begin{equation*}
  A_n:=\!\!\sum_{(k,\ell)\in T_{k_n,\ell _n}}
       \begin{pmatrix}
        X_{k-1,\ell}^{(r_n)}\varepsilon_{k,\ell}^{(r_n)} \\
        X_{k,\ell-1}^{(r_n)}\varepsilon_{k,\ell}^{(r_n)}
       \end{pmatrix},\qquad
  B_n:=\!\!\sum_{(k,\ell)\in T_{k_n,\ell _n}}
               \begin{pmatrix}
                 \big (X_{k-1,\ell}^{(r_n)}\big )^2 & X_{k-1,\ell}^{(r_n)}
                        X_{k,\ell-1}^{(r_n)} \\
                 X_{k-1,\ell}^{(r_n)}X_{k,\ell-1}^{(r_n)} & \big
                 (X_{k,\ell-1}^{(r_n)}\big )^2
               \end{pmatrix}.
 \end{equation*}
Concerning the asymptotic behaviour of the random vector \ $A_n$ \
and random matrix \ $B_n$ \ we can formulate the following two
propositions.

\begin{Pro}
   \label{pro:pro1} \
If \ $0<|\alpha| <1, \ |\beta|=1-|\alpha|$ \ and \eqref{eq:eq1.7} holds then
\begin{equation*}
 n^{-2}{r_n}^{-1/2}\big(|\gamma _{r_n}|+|\delta _{r_n}|\big )^{1/2}
B_n \qmean \big (32|\alpha||\beta |\big )^{-1/2}\Psi _{\alpha, \beta}
\qquad \text{as \ $n\to\infty$.}  
\end{equation*}

If \  $|\alpha|\in \{0,1\}, \ |\beta|=1-|\alpha|$ \ and \eqref{eq:eq1.8}
holds then
\begin{equation*}
 n^{-2}r_n^{-1}\big |\gamma _{r_n}^2-\delta _{r_n}^2\big |^{1/2}
B_n  \qmean  \Theta _{\alpha,\beta,\omega}
\end{equation*}
as \ $n\to\infty$, \ where
\begin{equation}
    \label{eq:eq1.9}
\omega:=\lim_{n\to\infty}\omega _{r_n}, \qquad \qquad \omega_{r_n}:=\alpha
\frac{\gamma_{r_n}}{\delta_{r_n}}+\beta 
  \frac{\delta_{r_n}}{\gamma_{r_n}}.
\end{equation}
\end{Pro}

\begin{Pro}
   \label{pro:pro2} \
If \  $0<|\alpha| <1, \ |\beta|=1-|\alpha|$ \ and \eqref{eq:eq1.7} holds then
\begin{equation*}
n^{-1}r_n^{-1/4}\big(|\gamma _{r_n}|+|\delta _{r_n}|\big )^{1/4}
A_n\distr {\mathcal N}_2\Big (0,\big (32|\alpha||\beta |\big
)^{-1/2}\Psi _{\alpha,\beta}\Big ) \qquad \text{as \ $n\to\infty$.}
\end{equation*}

If \ $|\alpha|\in \{0,1\}, \  |\beta|=1-|\alpha|$ \ and
\eqref{eq:eq1.8} holds then 
\begin{equation*}
 n^{-1}r_n^{-1/2}\big |\gamma _{r_n}^2-\delta _{r_n}^2\big
 |^{1/4}A_n\distr{\mathcal N}_2\bigg (0,\Theta
 _{\alpha,\beta,\omega}\bigg) \qquad \text{as \ $n\to\infty$.}
\end{equation*}
\end{Pro}

In case \ $|\alpha|\in \{0,1\}, \  |\beta|=1-|\alpha|$, \ and \
$|\omega|\ne 1$, \ $\Theta _{\alpha,\beta,\omega}$ \ is a regular
matrix, so Propositions \ref{pro:pro1} and \ref{pro:pro2} imply the
corresponding statement of Theorem \ref{main}. In the 
case \  $0<|\alpha| <1, \ |\beta|=1-|\alpha|$ \ we have \
$B_n^{-1}=\bar B_n/{\det B_n}$, \ and in this situation the statement of
Theorem \ref{main} is a consequence of the following propositions.  

\begin{Pro}
   \label{pro:pro3} \
If \  $0<|\alpha| <1, \ |\beta|=1-|\alpha|$ \ and \eqref{eq:eq1.7}
holds then
\begin{equation*}
n^{-4}r_n^{-1/2}\big(|\gamma _{r_n}|+|\delta _{r_n}|\big )^{1/2}
\det B_n \qmean 2\big(8|\alpha||\beta|\big)^{-3/2} \qquad
\text{as \ $n\to\infty$.} 
\end{equation*}
\end{Pro}

\begin{Pro}
   \label{pro:pro4} \
If \  $0<|\alpha| <1, \ |\beta|=1-|\alpha|$ \ and \eqref{eq:eq1.7}
holds then
\begin{equation*}
n^{-3}r_n^{-1/2}\big(|\gamma _{r_n}|+|\delta _{r_n}|\big )^{1/2}
\bar B_n A_n \distr {\mathcal N}_2\Big (0,\big (2(8 \alpha\beta)^2\big
)^{-1}\bar \Psi_{\alpha,\beta}\Big ) \qquad \text{as \ $n\to\infty$.}
\end{equation*}
\end{Pro}

Obviously, in the case \ $0\leq |\alpha|\leq 1, \ |\beta|=1-|\alpha|$
\ if \ $n$ \ is large enough, the corresponding sequences \
$\alpha_{r_n}$ \ and \ $\beta_{r_n}$ \ have the same signs as \
$\alpha$ \  and \ $\beta$, \ respectively. Hence, 
similarly to \cite{bpz3}, it suffices to prove Propositions
\ref{pro:pro3} and \ref{pro:pro4} for \ $0< \alpha,\beta<1,\ 
\alpha+\beta =1$.

\section{Covariance structure}
   \label{sec:sec2}
Let \ $\{X_{k,\ell}:k,\ell\in{\mathbb Z}\}$ \ be a stationary solution
of equation
 \eqref{eq:eq1.2} with parameters \ $(\alpha,\beta)$, \
 $|\alpha|+|\beta |<1$. \
Clearly \
$\cov(X_{i_1,j_1},X_{i_2,j_2})=\cov(X_{i_1-i_2,j_1-j_2},X_{0,0})$ \ for
 all \ $i_1,j_1,i_2,j_2\in{\mathbb Z}$. \
Let \ $R_{k,\ell}:=\cov(X_{k,\ell},X_{0,0})$ \ for \
$k,\ell\in{\mathbb Z}$. \
The following lemma is a natural generalization of Lemma 4 of
\cite{bpz2} (see also \cite{br2}).
\begin{Lem}
   \label{lem:lem1}
Let \ $\alpha\ne 0$ \ and \ $\beta \ne 0$. \ 
If \ $k,\ell\in{\mathbb Z}$ \ with \ $k\cdot\ell\leq0$ \ then
 \begin{equation}
    \label{eq:eq2.1}
  R_{k,\ell}
  =\sigma ^2_{\alpha,\beta}
    \left(\frac{1+\alpha ^2-\beta ^2-\sigma
        ^{-2}_{\alpha,\beta}}{2\alpha}\right)^{|k|}
    \left(\frac{2\beta}{1+\beta ^2-\alpha ^2 +\sigma
        ^{-2}_{\alpha,\beta}}\right)^{|\ell|}.
 \end{equation}
If \ $k,\ell\in{\mathbb Z}$ \ with \ $k\cdot \ell\geq 0$ \ then
 \begin{equation}
   \label{eq:eq2.2}
  R_{k,\ell}
  =R_{0,|k-\ell|}
   -\sum_{i=0}^{|k|\land|\ell|-1}\binom{|k-\ell|+2i}{i}\alpha^i
 \beta^{|k-\ell|+i}.
 \end{equation}
\end{Lem}

\begin{Rem}
   \label{rem:rem3}
If \ $\alpha >0$ \ and \ $\beta > 0$ \ then \ $R_{k,\ell}\geq 0$. \ If
\ $\alpha < 0$ \ or \ $\beta < 0$ \ we have
\begin{equation*}
0\leq |R_{k,\ell}|\leq \widetilde R_{k,\ell}:=\cov (\widetilde
X_{k,\ell}, \widetilde X_{0,0}), \qquad k,\ell\in{\mathbb Z},
\end{equation*}
where \ $\{\widetilde X_{k,\ell}:k,\ell\in{\mathbb Z}\}$ \ is a
stationary solution of equation \eqref{eq:eq1.2} with parameters \
$(|\alpha|,|\beta|)$. \
\end{Rem}

Besides representations \eqref{eq:eq2.1} and \eqref{eq:eq2.2} one can
express the covariances as special cases of Appell's hypergeometric
series \ $F_4(a,b,c,d;x,y)$ \ defined by
\begin{equation*}
F_4(a,b,c,d;x,y):=\sum_{m=0}^{\infty }\sum_{n=0}^{\infty }
\frac {(a)_{m+n}(b)_{m+n}}{(c)_m(d)_nm!n!}x^my^n, \qquad
\sqrt {|x|}+\sqrt {|y|}< 1,
\end{equation*}
where \ $a,b,c,d\in {\mathbb N}$ \ and \ $(a)_n:=a(a+1)\dots (a+n-1)$ \
\cite{be}.
\begin{Lem}
   \label{lem:lem2}
Let \ $\alpha\ne 0$ \ and \ $\beta \ne 0$. \ 
If \ $k,\ell\in{\mathbb Z}$ \ with \ $k\cdot\ell\leq0$ \ then
 \begin{equation}
    \label{eq:eq2.3}
R_{k,\ell}=\alpha ^{|k|}\beta ^{|\ell |}F_4
\big(|k|+1, |\ell |+1,|k|+1,|\ell |+1; \alpha ^2,\beta ^2\big ).
\end{equation}
If \ $k,\ell\in{\mathbb Z}$ \ with \ $k\cdot \ell\geq 0$ \ then
 \begin{equation*}
R_{k,\ell}=\alpha ^{|k|}\beta ^{|\ell |}\binom{|k|+|\ell |}{|k|}F_4
\big(|k|+|\ell |+1,1,|k|+1,|\ell |+1; \alpha ^2,\beta ^2\big ).
\end{equation*}
Moreover, in this case we have
 \begin{equation}
   \label{eq:eq2.4}
R_{k,\ell}=\big(\sign(\alpha)\big)^{|k|}\big(\sign(\beta)\big)^{|\ell|}
\sum _{i=0}^{\infty}\big (|\alpha |+|\beta|\big )^{|k|+|\ell |+2i}
{\mathsf P}\big (S^{(\nu)}_{i,|k|+|\ell|+i}=|\ell |+i\big ),
  \end{equation}
where \ $S^{(\nu)}_{n,m}:=S^{(\nu)}_n+S^{(1-\nu)}_m$, \
$\nu:=|\alpha|/\big(|\alpha|+|\beta|\big)$ \ and \ $S^{(\nu)}_n$ \ and
\ $S^{(1-\nu)}_m$ \ are independent binomial random variables with
parameters \ $(n,\nu)$ \ and \ $(m,1-\nu)$, respectively.
\end{Lem}

\bigskip
\noindent
{\bf Proof.} \
 The statements directly follow from representation
\eqref{eq:eq1.3} and from the 
independence of the error terms \ $\varepsilon _{i,j}$. \proofend

\bigskip
We remark, that as
\begin{equation*}
F_4\Big (a,b,a,b;\frac {-x}{(1-x)(1-y)},\frac {-y}{(1-x)(1-y)}\Big )=
\frac{(1-x)^b(1-y)^a}{1-xy},
\end{equation*}
representation \eqref{eq:eq2.1} directly follows from \eqref{eq:eq2.3}.

\begin{Pro}
  \label{diffbound}
If \ $\alpha\beta>0$, \ $|\alpha|+|\beta|<1$ \ then there exists a
universal positive constant \ $K$ \ 
such that
\begin{equation*}
\big|R_{k-1,\ell +1}-R_{k,\ell}\big |\leq \frac
K{(\alpha\beta)^{3/2}}, \qquad \qquad k,\ell\in{\mathbb Z}.
\end{equation*}
\end{Pro}

\bigskip
\noindent
{\bf Proof.} \ Without loss of generality we may assume \ $\alpha>0$ \
and \ $\beta>0$. \

Suppose \ $k>0, \ \ell\geq 0$, \ so \ $(k-1)(\ell
+1)\geq 0$ \ and \ $k\cdot \ell\geq 0$. \ Using notations introduced
in Lemma \ref{lem:lem2} with the help of \eqref{eq:eq2.4} we obtain
\begin{equation}
   \label{eq:eq2.5}
R_{k-1,\ell +1}-R_{k,\ell}=\sum
_{i=0}^{\infty}(\alpha+\beta)^{k+\ell+2i}
\Delta_{k,\ell,i}(\nu),
\end{equation}
where
\begin{equation*}
\Delta_{i,k,\ell}(\nu):={\mathsf P}\big (S^{(\nu)}_{i,k+\ell+i}=\ell
+i+1\big )-{\mathsf P}\big (S^{(\nu)}_{i,k+\ell+i}=\ell +i\big ).
\end{equation*}
According to Theorem 2.6 of \cite{bpz3} \ $\Delta_{i,k,\ell}(\nu)$ \
can be approximated by
\begin{align*}
\widetilde\Delta_{i,k,\ell}(\nu):=\frac
1{\big(2\pi\nu(1-\nu)(k+\ell+2i)\big)^{1/2}}\Bigg (\exp &\bigg\{-\frac
  {\big(\nu\ell -(1-\nu)k+1\big)^2}{2\nu(1-\nu)(k+\ell+2i)}\bigg\} \\
&-\exp \bigg\{-\frac
  {\big(\nu\ell -(1-\nu)k\big)^2}{2\nu(1-\nu)(k+\ell+2i)}\bigg\}\Bigg)
\end{align*}
where
\begin{equation*}
\big| \widetilde \Delta_{i,k,\ell}(\nu)-\Delta_{i,k,\ell}(\nu)\big|
\leq \frac{\widetilde C}{\big (\nu(1-\nu)(k+\ell+2i)\big )^{3/2}}
\end{equation*}
with some positive constant \ $\widetilde C$. \ Thus, if in the right hand side
of \eqref{eq:eq2.5} we replace \ $\Delta_{i,k,\ell}(\nu)$ \
with \ $\widetilde\Delta_{i,k,\ell}(\nu)$, \ the error of the approximation is
\begin{equation*}
\sum _{i=0}^{\infty}(\alpha+\beta)^{k+\ell+2i}
\big| \widetilde \Delta_{i,k,\ell}(\nu)-\Delta_{i,k,\ell}(\nu)\big|
\leq  \frac {\widetilde C}{\big (\nu(1-\nu)\big )^{3/2}}\zeta(3/2)\leq
 \frac C{(\alpha\beta)^{3/2}},
\end{equation*}
where \ $\zeta(x)$ \ denotes Riemann's zeta function.

To find an upper bound for the approximating sum consider first the
case \ $\nu\ell -(1-\nu)k\geq 0$. \ In this case
\begin{align*}
\sum _{i=0}^{\infty}(\alpha&\!+\!\beta)^{k+\ell+2i}
\big| \widetilde \Delta_{i,k,\ell}(\nu)\big | \!\leq\! \sum
_{i=0}^{\infty}\! \frac{2(\nu\ell
  \!-\!(1\!-\!\nu)k)\!+1\!}{\pi
  ^{1/2}\big(2\nu(1\!-\!\nu)(k\!+\!\ell\!+\!2i)\big )^{3/2}}
\exp \bigg\{\!\! -\!\frac
  {\big(\nu\ell
    \!-\!(1\!-\!\nu)k\big)^2}{2\nu(1\!-\!\nu)(k\!+\!\ell\!+\!2i)}\!\bigg\} \\
&\leq \frac {\zeta (3/2)+1}{\big
  (\nu(1\!-\!\nu)\big )^{3/2}} \!+\!\frac 1{2\nu(1\!-\!\nu)}\widetilde \Phi
\Bigg ( \frac {\nu\ell -(1-\nu)k}{\big
  (2\nu(1\!-\!\nu)(k\!+\!\ell)\big )^{1/2}}\Bigg )
\!\leq\! \frac {\zeta (3/2)+2}{\big  (\nu(1\!-\!\nu)\big )^{3/2}}\!\leq\!
\frac {\zeta (3/2)+2}{(\alpha \beta )^{3/2}},
\end{align*}
where $\widetilde\Phi (x)$ is the error function defined by
\begin{equation*}
\widetilde\Phi(x):=\frac 2{\pi ^{1/2}}\int\limits_0^x {\mathrm
  e}^{-t^2/2}{\mathrm d}t, \qquad x>0.
\end{equation*}
Case \ $\nu\ell -(1-\nu)k<0$ \ follows by symmetry.

In case \ $k\leq0, \ \ell <0$ \ implying \ $(k-1)(\ell +1)\geq 0$ \
and \ $k\cdot \ell >0$, \ we have
\begin{equation*}
R_{k-1,\ell +1}-R_{k,\ell}=\sum
_{i=0}^{\infty}(\alpha+\beta)^{-k-\ell+2i}\Big ( {\mathsf P}\big
(S^{(\nu)}_{i,-k-\ell+i}=-\ell +i-1\big )-{\mathsf P}\big
(S^{(\nu)}_{i,-k-\ell+i}=-\ell +i\big )\Big )
\end{equation*}
and the statement can be proved similarly to the previous case.

Now, suppose \ $k>0, \ \ell <0$, \ so \ $(k-1)(\ell +1)\leq 0$ \ and
\ $k\cdot \ell \leq 0$. \ Using the form \eqref{eq:eq2.1} of the
covariances direct calculations show 
\begin{equation*}
R_{k-1,\ell +1}-R_{k,\ell}=R_{k,\ell }\frac {1-(\alpha+\beta)^2+\sigma
        ^{-2}_{\alpha,\beta}}{2\alpha\beta}.
\end{equation*}
It is not difficult to see that \ $1-(\alpha +\beta )^2\leq \sigma
^{-2}_{\alpha,\beta}$, \ so we have
\begin{equation*}
\big |R_{k-1,\ell +1}-R_{k,\ell}\big |\leq \big |R_{k,\ell}\big |\,
\frac {\sigma ^{-2}_{\alpha,\beta}}{\alpha\beta}\leq \frac
1{\alpha\beta}.
\end{equation*}
In a similar way one can obtain the result for \ $k\leq 0, \ \ell
\geq 0$ \ that completes the proof. \proofend

\bigskip
Using the notations of Lemma \ref{lem:lem2} with the help of the
exponential approximation one can easily have the analogue of 
Corollary 2.7 of \cite{bpz3}.

\begin{Cor}
   \label{cor:cor1} 
If \ $\alpha\beta>0$, \ $|\alpha|+|\beta|<1$ \ then there exists a
constant \ $C>0$  \ such that for all \ $k,\ell>1$ \ and \ $0\leq
i\leq k+\ell -1$ we have
\begin{equation*}
\Big | {\mathsf P}\big (S^{(\nu)}_{k,\ell}=i+1\big )-{\mathsf P}\big
(S^{(\nu)}_{k,\ell}=i\big ) \Big |\leq \frac C{\alpha\beta(k+\ell)}.
\end{equation*}
\end{Cor}

\begin{Rem}
   \label{rem:rem4} 
Using Theorem 2.4 of \cite{bpz3} it is not difficult to show that under
conditions of Corollary \ref{cor:cor1} there exists a
constant \ $D>0$  \ such that for all \ $k,\ell>1$ \ and \ $0\leq
i\leq k+\ell$ we have
\begin{equation*}
\Big | {\mathsf P}\big (S^{(\nu)}_{k,\ell}=i\big )\Big |\leq 
\frac D{\alpha\beta(k+\ell)^{1/2}}.
\end{equation*}
\end{Rem}

Now, let \ $\{X^{(n)}_{k,\ell}:k,\ell\in{\mathbb Z}\}, \ n\in{\mathbb
  N}$, \ be a nearly unstable sequence of stationary processes described
in Theorem \ref{main}. For each \ $n\in{\mathbb N}$ \ let us introduce
the piecewise constant random fields
\begin{alignat*}{2}
Z_{1,0}^{(n)}(s,t)&:=r_n^{-1/4}X_{[ns]+1,[nt]}^{(r_n)}, \qquad
Z_{0,1}^{(n)}(s,t)&:=r_n^{-1/4}X_{[ns],[nt]+1}^{(r_n)}, \quad
\phantom{s,t\in{\mathbb R}.} \\
Y_{1,0}^{(n)}(s,t)&:=r_n^{-1/2}X_{[ns]+1,[nt]}^{(r_n)}, \qquad
Y_{0,1}^{(n)}(s,t)&:=r_n^{-1/2}X_{[ns],[nt]+1}^{(r_n)}, \quad
s,t\in{\mathbb R}.
\end{alignat*}

\begin{Pro}
  \label{covlim}
Let \ $s_1,t_1,s_2,t_2\in{\mathbb R}$.

If \ $0<|\alpha |<1, \ |\beta|=1-|\alpha|$ \ and \eqref{eq:eq1.7}
holds then for all \
$(i_1,j_1),(i_2,j_2)\in \big\{ (1,0),(0,1)\big\}$ \ we have
 \begin{alignat*}{2}
  &\lim_{n\to\infty}
   \big(|\gamma_{r_n}|\!+\!|\delta _{r_n}|\big )^{1/2}\cov
  \big (Z_{i_1,j_1}^{(n)}(s_1,t_1),Z_{i_2,j_2}^{(n)}(s_2,t_2)\big)
  =0\ \
  &\text{if  $s_1\!-\!s_2\!\ne\! t_1\!-\!t_2$,}\\[2mm]
  &\limsup_{n\to\infty}
     \big(|\gamma_{r_n}|\!+\!|\delta _{r_n}|\big )^{1/2}
    \Big |\cov \big
    (Z_{i_1,j_1}^{(n)}(s_1,t_1),Z_{i_2,j_2}^{(n)}(s_2,t_2)\big )\Big|
   \!\leq\!\frac 1{\sqrt {8|\alpha||\beta|}}\ \
  &\text{if  $s_1\!-\!s_2\!=\!t_1\!-\!t_2$.}
 \end{alignat*}

If \ $|\alpha |\in \{0,1\}, \ |\beta|=1-|\alpha|$ \ and
\eqref{eq:eq1.8} holds then for all \
$(i_1,j_1),(i_2,j_2)\!\in\! \big\{ (1,0),(0,1)\big\}$ \ we have
 \begin{alignat*}{2}
  &\lim_{n\to\infty}
    \big |\gamma_{r_n}^2-\delta _{r_n}^2\big |^{1/2}\cov
  \big (Y_{i_1,j_1}^{(n)}(s_1,t_1),Y_{i_2,j_2}^{(n)}(s_2,t_2)\big)
  =0\quad
  &\text{if \ $s_1\!-\!s_2\ne t_1\!-\!t_2$,}\\[2mm]
  &\limsup_{n\to\infty}
    \big|\gamma_{r_n}^2-\delta _{r_n}^2\big |^{1/2}
    \Big |\cov \big
    (Y_{i_1,j_1}^{(n)}(s_1,t_1),Y_{i_2,j_2}^{(n)}(s_2,t_2)\big )\Big|
   \leq\frac 12\quad
  &\text{if \ $s_1\!-\!s_2=t_1\!-\!t_2$.}
 \end{alignat*}

Moreover, if \ $s_1-s_2\ne t_1-t_2$ \ then the 
convergence to \ $0$ \ in both cases has an exponential rate.
\end{Pro}

\bigskip
\noindent
{\bf Proof.} \ For simplicity we consider only the case \ $0\leq
\alpha,\beta \leq 1$. \ The other cases can be handled in a similar way.

First, let \ $0<\alpha <1$, \ so \ $\beta =1-\alpha$. \ Without loss of
generality we may assume \ $\alpha _{r_n} > 0$, \ $\beta _{r_n} >
0$ \ and \ $\delta _{r_n} > 0$, \ $\gamma _{r_n} >0$. \ As
\begin{equation*}
r_n^{-1/2} \sigma ^2_{\alpha _{r_n},\beta _{r_n}}\!=\bigg (\big (\gamma
_{r_n}+\delta _{r_n}\big)\Big (2-\frac {\gamma_{r_n}+\delta
  _{r_n}}{r_n}\Big)\Big (2\alpha -\frac {\gamma_{r_n}-\delta
  _{r_n}}{r_n}\Big)\Big(2(1-\alpha)+\frac {\gamma_{r_n}-\delta
  _{r_n}}{r_n}\Big) \bigg )^{-1/2}
\end{equation*}
we have
\begin{equation}
    \label{eq:eq2.6}
 \lim_{n\to\infty}\big(\gamma_{r_n}+\delta
 _{r_n}\big)^{1/2}r_n^{-1/2}\sigma ^2_{\alpha _{r_n},\beta _{r_n}}
 =\frac 1{\sqrt{8\alpha(1-\alpha)}}=\frac 1{\sqrt{8\alpha\beta}}.
\end{equation}

Suppose \ $s_1-s_2\geq 0\geq t_1-t_2$, \ so \ $[ns_1]-[ns_2]\geq 0\geq
[nt_1]-[nt_2]$. \ By \eqref{eq:eq2.1}
\begin{equation*}
0\leq \cov
  \big (Z_{1,0}^{(n)}(s_1,t_1),Z_{1,0}^{(n)}(s_2,t_2)\big)
  \leq r_n^{-1/2}
  \sigma ^2_{\alpha _{r_n},\beta _{r_n}} \Big (1-\frac 1{\varrho
    _{r_n}}\Big )^{\frac n2 |s_1-s_2|}\Big (1+\frac 1{\tau
    _{r_n}}\Big )^{-\frac n2|t_1-t_2|}
\end{equation*}
if \ $n$ \ is large enough, where
\begin{equation}
  \label{eq:eq2.7}
\varrho _{r_n}\!:=\!\frac {2\alpha _{r_n}}{2\alpha _{r_n}\!-\!1\!-\!\alpha
  _{r_n}^2\!+\!\beta _{r_n}^2\!+\!\sigma ^{-2}_{\alpha _{r_n},\beta _{r_n}}},
\qquad
\tau _{r_n}\!:=\!\frac {2\beta _{r_n}}{1\!+\!\beta_{r_n}^2\!-\!\alpha
  _{r_n}^2+\sigma ^{-2}_{\alpha _{r_n},\beta _{r_n}}-2\beta _{r_n}}.
\end{equation}
As
\begin{equation*}
\sigma ^2_{\alpha,\beta}=\big((1+\alpha^2-\beta^2)^2-4\alpha ^2\big)^{-1/2},
\end{equation*}
it is easy to see that \ $\varrho _{r_n}\to \infty$ \ and \ $\tau _{r_n}\to
\infty$ \ as \ $n\to\infty$. \ Moreover, condition \eqref{eq:eq1.7} ensures
that \ $n\varrho _{r_n}^{-1}\to\infty$ \ and \ $n\tau
_{r_n}^{-1}\to\infty$ \ as \
$n\to\infty$. \ Hence, if \ $s_1=s_2$ \ and \ $t_1=t_2$, \
\begin{equation*}
\lim _{n\to\infty}\big(\gamma_{r_n}+\delta
 _{r_n}\big)^{1/2}\cov
  \big (Z_{1,0}^{(n)}(s_1,t_1),Z_{1,0}^{(n)}(s_2,t_2)\big)=\frac
 1{\sqrt{8\alpha\beta}},  
\end{equation*}
otherwise it converges to \ $0$ \ in exponential rate.

Further, let \ $s_1-s_2>0$ \ and \ $t_1-t_2>0$. \ In this case \
$[ns_1]-[ns_2]\geq 0$ \ and \ $[nt_1]-[nt_2]\geq 0$, \ so by
\eqref{eq:eq2.2} we have
\begin{equation}
    \label{eq:eq2.8}
0\!\leq\! \cov
  \big (Z_{1,0}^{(n)}(s_1,t_1),Z_{1,0}^{(n)}(s_2,t_2)\big)\!\leq\! r_n^{-1/2}
  \sigma ^2_{\alpha _{r_n},\beta _{r_n}} \Big (1\!+\!\frac 1{\tau
    _{r_n}}\Big )^{-\left| [ns_1]-[ns_2]-[nt_1]+[nt_2]\right|}.
\end{equation}
If \ $s_1-s_2\ne t_1-t_2$ \ then similarly to the previous case one can
show that the right hand side of \eqref{eq:eq2.8} converges to \ $0$ \
in exponential rate as \ $n\to\infty$.

In case  \ $s_1-s_2=t_1-t_2$ \ we have \ $\big
|[ns_1]-[ns_2]-[nt_1]+[nt_2]\big |\leq 2$, \ so by  \eqref{eq:eq2.8}
\begin{equation*}
\limsup _{n\to\infty}\big(\gamma_{r_n}+\delta
 _{r_n}\big)^{1/2}\cov
  \big (Z_{1,0}^{(n)}(s_1,t_1),Z_{1,0}^{(n)}(s_2,t_2)\big)\leq\frac
 1{\sqrt{8\alpha\beta}}. 
\end{equation*}
Obviously, the same results hold for the covariances \ 
$\cov\big (Z_{1,0}^{(n)}(s_1,t_1),Z_{0,1}^{(n)}(s_2,t_2)\big)$, \break
$\cov\big (Z_{0,1}^{(n)}(s_1,t_1),Z_{1,0}^{(n)}(s_2,t_2)\big)$ \ and \
$\cov\big (Z_{0,1}^{(n)}(s_1,t_1),Z_{0,1}^{(n)}(s_2,t_2)\big)$.

\smallskip
Now, consider for example the case \ $\alpha=1, \ \beta=0$. \ Without loss
of generality we may assume \ $\alpha _{r_n}> 0$. \ Furthermore, \ $|\alpha
_{r_n}|+|\beta _{r_n}|<1$ \ implies  \ $\gamma_{r_n}> 0$ \ and \
$|\delta _{r_n}|<\gamma _{r_n}$. \ As
\begin{equation*}
r_n^{-1} \sigma ^2_{\alpha _{r_n},\beta _{r_n}}=\bigg (\big( \gamma
_{r_n}^2-\delta _{r_n}^2\big)\Big (2-\frac {\gamma_{r_n}+\delta
  _{r_n}}{r_n}\Big)\Big (2 -\frac {\gamma_{r_n}-\delta
  _{r_n}}{r_n}\Big) \bigg )^{-1/2}
\end{equation*}
we have
\begin{equation}
    \label{eq:eq2.9}
 \lim_{n\to\infty}\big(\gamma_{r_n}^2-\delta
 _{r_n}^2\big)^{1/2}r_n^{-1}\sigma ^2_{\alpha _{r_n},\beta _{r_n}}
 =\frac 12.
\end{equation}

Again, suppose \ $s_1-s_2\geq 0\geq t_1-t_2$. \ The form of covariances
\eqref{eq:eq2.1} implies that if \ $n$ \ is large enough
\begin{equation}
   \label{eq:eq2.10}
0\!\leq\! \Big |\cov
  \big (Y_{1,0}^{(n)}(s_1,t_1),Y_{1,0}^{(n)}(s_2,t_2)\big)\Big
  |\!\leq\! r_n^{-1} 
  \sigma ^2_{\alpha _{r_n},\beta _{r_n}} \Big (1\!-\!\frac 1{\varrho
    _{r_n}}\Big )^{\frac n2 |s_1-s_2|}\Big (1\!+\!\frac 1{|\tau
    _{r_n}|}\Big )^{-\frac n2 |t_1-t_2|},
\end{equation}
where \ $\varrho _{r_n}$ \ and \ $\tau _{r_n}$ \ are defined by
\eqref{eq:eq2.7}.
Obviously, if \ $s_1=s_2$ \ and \ $t_1=t_2$ \ then \eqref{eq:eq2.9} implies
\begin{equation}
   \label{eq2.11}
\limsup _{n\to\infty}\big(\gamma_{r_n}^2-\delta
_{r_n}^2\big)^{1/2}\Big |\cov
  \big (Y_{1,0}^{(n)}(s_1,t_1),Y_{1,0}^{(n)}(s_2,t_2)\big)\Big |\leq\frac 12.
\end{equation}
Further, we have \ $\varrho _{r_n}\to\infty $ \ as \
$n\to\infty$ \ and now \eqref{eq:eq1.8} ensures \
$n\varrho_{r_n}^{-1}\to\infty$. \ Thus, as \ $1+1/|\tau _{r_n}|\geq
1$, \ if \ $s_1\ne s_2$ \ then
\begin{equation}
   \label{eq:eq2.12}
\big(\gamma_{r_n}^2-\delta _{r_n}^2\big)^{1/2}\Big |\cov
  \big (Y_{1,0}^{(n)}(s_1,t_1),Y_{1,0}^{(n)}(s_2,t_2)\big)\Big |\to 0
\end{equation}
as \ $n\to\infty$ \ in exponential rate. Now, let us assume \
$s_1=s_2$ \ and \ $t_1\ne t_2$. \ Short calculation shows 
\begin{equation}
   \label{eq:eq2.13}
\Big (1+\frac 1{|\tau _{r_n}|}\Big )^{-1}=\frac {2|\delta
  _{r_n}|}{2\gamma _{r_n}-\frac {\gamma _{r_n}^2-\delta
    _{r_n}^2}{r_n}+
\big (\gamma _{r_n}^2-\delta _{r_n}^2\big )^{1/2}\Big (\frac {\gamma
  _{r_n}^2-\delta _{r_n}^2}{r_n} -4\frac {\gamma_{r_n}}{r_n}+4 \Big )^{1/2}}.
\end{equation}
If \ $|\delta| <\gamma$ \ then
\begin{equation*}
\lim_{n\to\infty}\Big (1+\frac 1{|\tau _{r_n}|}\Big )^{-1}=
\frac{|\delta|}{\gamma +(\gamma^2-\delta ^2)^{1/2}}<1,
\end{equation*}
so using  \eqref{eq:eq2.9} and \eqref{eq:eq2.10} we obtain again
\eqref{eq:eq2.12}.
Further, condition \eqref{eq:eq1.8} implies
\begin{equation*}
\lim_{n\to\infty}n\big (\gamma _{r_n}^2-\delta _{r_n}^2\big
)^{1/2}=\infty .
\end{equation*}
Hence, with the help of \eqref{eq:eq2.13} one can easily see that if
\ $|\delta |=\gamma\ne 0$, \ or \ $\delta =\gamma=0$ \ and \ $\lim
_{n\to\infty} \gamma_{r_n}|\delta _{r_n}|^{-1}=1$, \ we obtain \
$|\tau_{r_n}| \to\infty $ \ and \ $n|\tau_{r_n}|^{-1}\to\infty$ \ as
\ $n\to \infty$. \ Thus, \eqref{eq:eq2.9} and \eqref{eq:eq2.10} imply
\eqref{eq:eq2.12} and the rate of convergence is again exponential.
In case \  $\delta =\gamma=0$ \ and \ $\lim _{n\to\infty}
\gamma_{r_n}|\delta _{r_n}|^{-1}=|\omega| >1$ \ we have
\begin{equation*}
\lim_{n\to\infty}\Big (1+\frac 1{|\tau _{r_n}|}\Big )^{-1}=
\frac{1}{|\omega| +(\omega^2-1)^{1/2}}<1,
\end{equation*}
that implies  \eqref{eq:eq2.12}. Finally, if \ $\delta =\gamma=0$ \
and \ $\lim _{n\to\infty} \gamma_{r_n}|\delta _{r_n}|^{-1}=\infty$ \ then
\eqref{eq:eq2.12} follows from
\begin{equation*}
\lim_{n\to\infty}\Big (1+\frac 1{|\tau _{r_n}|}\Big )^{-1}=0.
\end{equation*}

Now, let \ $s_1-s_2>0$ \ and \ $t_1-t_2>0$. \ Lemma \ref{lem:lem1} and
Remark \ref{rem:rem3} imply
\begin{equation*}
0\leq \Big |\cov
  \big (Y_{1,0}^{(n)}(s_1,t_1),Y_{1,0}^{(n)}(s_2,t_2)\big)\Big |\leq r_n^{-1}
  \sigma ^2_{\alpha _{r_n},\beta _{r_n}} \Big (1+\frac 1{|\tau
    _{r_n}|}\Big )^{-\left |[ns_1]-[ns_2]-[nt_1]+[nt_2]\right |},
\end{equation*}
where \ $\tau_{r_n}$ \ is defined by \eqref{eq:eq2.7}. If \
$s_1-s_2=t_1-t_2$ \ then as \ $\big |[ns_1]-[ns_2]-[nt_1]+[nt_2]\big |\leq
2$ \ and \ $1+1/|\tau _{r_n}|\geq 1$, \ using \eqref{eq:eq2.9} we obtain
\eqref{eq2.11}. Finally, if \ $s_1-s_2\ne t_1-t_2$ \ then to prove
\eqref{eq2.11} one has to do the same considerations as in the case \
$s_1=s_2$ \ and \ $t_1\ne t_2$. \ 

\space \hfill\proofend

In order to estimate the covariances we make use of the following
lemma which is a natural generalization of Lemma 2.8 of \cite{bpz3}.

\begin{Lem}
   \label{lem:lem3}
Let \ $\xi_1, \xi_2, \ldots $ \ be independent random variables with \
${\mathsf E}\xi_i=0$, \ ${\mathsf E}\xi_i^2=1$ \ for all \ $i\in
{\mathbb N}$, \ and \ $M_4:=\sup_{i\in {\mathbb N}}{\mathsf E}\xi_i^4<\infty$.
\ Let \ $a_1,a_2,\ldots ,b_1, b_2,\ldots,c_1,c_2\ldots$,
$d_1,d_2\ldots \in{\mathbb R}$, \ such that \
$\sum_{i=1}^{\infty}a_i^2<\infty, \ \sum_{i=1}^{\infty}b_i^2<\infty, \
\sum_{i=1}^{\infty}c_i^2<\infty$ \ and \
$\sum_{i=1}^{\infty}d_i^2<\infty$. Let
$$X:=\sum_{i=1}^{\infty}a_i\xi_i, \quad Y:=\sum_{i=1}^{\infty}b_i\xi_i,
 \quad Z:=\sum_{i=1}^{\infty}c_i\xi_i, \quad W:=\sum_{i=1}^{\infty}d_i\xi_i,$$
where the convergence of the infinite sums is understood in \ ${\mathsf
  L}_2$-sense.
Then
\begin{equation}
   \label{eq:eq2.14}
\cov(XY,ZW)=\sum_{i=1}^{\infty} ({\mathsf E}
\xi _i^4-3)\,  a_ib_ic_id_i+\cov (X,Z)\cov(Y,W)+\cov (X,W)\cov(Y,Z).
\end{equation}
Moreover, if  \ $a_i,b_i,c_i,d_i\geq 0$ \ then
 $$0\leq\cov (XY,ZW)\leq M_4\cov (X,Z)\cov(Y,W)+M_4\cov (X,W)\cov(Y,Z),$$
and
$$0\leq {\mathsf E} XYZW\leq M_4\big({\mathsf E} XZ\,{\mathsf E}
YW+{\mathsf E} XW\,{\mathsf E} YZ +{\mathsf E} XY\,{\mathsf E} ZW \big).$$
\end{Lem}

\begin{Rem}
  \label{rem:rem5}
Using the definitions of Lemma \ref{lem:lem3} from \eqref{eq:eq2.14}
one can easily see, that
\begin{equation*}
\big |\cov(XY,ZW) \big| \leq \cov(\widetilde X\widetilde Y,\widetilde
Z\widetilde W),
\end{equation*}
where
\begin{equation*}
\widetilde X:=\sum_{i=1}^{\infty}|a_i|\xi_i, \quad \widetilde
Y:=\sum_{i=1}^{\infty} |b_i|\xi_i,
 \quad \widetilde Z:=\sum_{i=1}^{\infty}|c_i|\xi_i, \quad \widetilde
 W:=\sum_{i=1}^{\infty}|d_i|\xi_i.
\end{equation*}
\end{Rem}

\section{Proof of Proposition \ref{pro:pro1}}
   \label{sec:sec3}
Let us assume \ $\alpha _{r_n}\ne 0$ \ and \ $\beta _{r_n}\ne 0$.
\ Using the stationarity of \ $\big\{ X_{k,\ell}^{(r_n)}: k,\ell \in
{\mathbb Z}\big\}$ \ and Lemma \ref{lem:lem1} we obtain
\begin{align*}
{\mathsf E}B_n&=\sum_{(k,\ell)\in T_{k_n,\ell _n}}
               \begin{pmatrix}
                 \var \big (X_{0,0}^{(r_n)}\big ) & \cov
                 \big(X_{0,0}^{(r_n)},
                        X_{1,-1}^{(r_n)}\big ) \\
                \cov \big (X_{0,0}^{(r_n)},X_{1,-1}^{(r_n)}\big ) &
                \var \big (X_{0,0}^{(r_n)}\big )
               \end{pmatrix} \\
=&\frac{(k_n+\ell _n)(k_n+\ell _n+1)}2\sigma ^2_{\alpha _{r_n},\beta
  _{r_n}} \begin{pmatrix}
              1&  D_{r_n}\\ D_{r_n} & 1
          \end{pmatrix}=\frac{n(n+1)}2\sigma ^2_{\alpha _{r_n},\beta
  _{r_n}} \begin{pmatrix}
              1&  D_{r_n}\\ D_{r_n} & 1
          \end{pmatrix},
\end{align*}
where
\begin{equation*}
D_{r_n}= \left(\frac{1+\alpha_{r_n}^2-\beta_{r_n}^2-\sigma
        ^{-2}_{\alpha_{r_n},\beta_{r_n}}}{2\alpha_{r_n}}\right)
    \left(\frac{2\beta_{r_n}}{1+\beta_{r_n}^2-\alpha_{r_n}^2 +\sigma
        ^{-2}_{\alpha_{r_n},\beta_{r_n}}}\right) .
\end{equation*}

If \ $0<|\alpha |<1$ \ and \ $|\beta |=1-|\alpha |$ \ then it is not difficult
to see that \ $\sigma
        ^{-2}_{\alpha_{r_n},\beta_{r_n}}\to 0$ \ and in this way \ $D_{r_n}\to
\sign(\alpha\beta) $ \ as \ $n\to\infty$. \ Hence, using the same arguments
as in the proof of \eqref{eq:eq2.6} we obtain
\begin{equation}
   \label{eq:eq3.1}
\lim _{n\to\infty} n^{-2}{r_n}^{-1/2}\big(|\gamma _{r_n}|+|\delta
_{r_n}|\big ) ^{1/2}
{\mathsf E}B_n  = \big (32|\alpha||\beta |\big )^{-1/2}\Psi _{\alpha, \beta}.
\end{equation}

If \ $|\alpha |\in \{0,1\}$ \ and \ $|\beta |=1-|\alpha |$, \ again, we have
\ $\sigma^{-2}_{\alpha_{r_n},\beta_{r_n}}\to 0$ \ as \ $n\to\infty$, \  and
similarly to the proof of \eqref{eq:eq2.9} one can see
\begin{equation*}
\lim _{n\to\infty} n^{-2}{r_n}^{-1}\big|\gamma _{r_n}^2-\delta
_{r_n}^2\big |^{1/2}\frac{n(n+1)}2\sigma^2_{\alpha_{r_n},\beta_{r_n}}=\frac14.
\end{equation*}
Concerning the limit of \ $D_{r_n}$ \ from the four possible cases
that can be handled in the same way we
consider only the case \ $\alpha =1, \ \beta =0$.  \ In this case \
$\alpha \frac {\gamma_{r_n}}{\delta _{r_n}}+\beta \frac
{\delta_{r_n}}{\gamma_{r_n}}=
\frac {\gamma_{r_n}}{\delta _{r_n}}$ \ and
we may assume \
$\alpha _{r_n}>0$ \ and thus \ $|\delta _{r_n}|\leq \gamma _{r_n}$ \
(hence \ $\gamma _{r_n}>0$). \
Obviously,
\begin{equation*}
\lim _{n\to\infty}\frac{1+\alpha_{r_n}^2-\beta_{r_n}^2-\sigma
        ^{-2}_{\alpha_{r_n},\beta_{r_n}}}{2\alpha_{r_n}}=1,
\end{equation*}
and
\begin{align*}
\frac{2\beta_{r_n}}{1+\beta_{r_n}^2-\alpha_{r_n}^2 +\sigma
        ^{-2}_{\alpha_{r_n},\beta_{r_n}}}=\Bigg (\frac {\gamma
        _{r_n}-\delta _{r_n}}{2r_n}&-\sign(\omega)\Big(1-\frac {\gamma
        _{r_n}-\delta _{r_n}}{2r_n}\Big)^{1/2}\\ \times\bigg (\frac
      {\gamma_{r_n}}{|\delta _{r_n}|}\Big(1-\frac {\gamma
        _{r_n}-\delta _{r_n}}{2r_n}\Big)^{1/2} 
&+\Big (\frac
      {\gamma ^2_{r_n}}{\delta ^2_{r_n}}-1\Big )^{1/2}\Big(1-\frac {\gamma
        _{r_n}+\delta _{r_n}}{2r_n}\Big)^{1/2}\bigg )\Bigg )^{-1}. 
\end{align*}
Hence,
\begin{equation*}
\lim_{n\to\infty} D_{r_n}=\begin{cases}
                          -\sign(\omega)\big(|\omega|+(\omega
                          ^2-1)^{1/2}\big )^{-1}& 
                          \text {if \ $|\omega |<\infty $,}\\
                          0& \text{if \ $|\omega|=\infty$,}
                          \end{cases}
\end{equation*}
where \ $\omega $ is the limit  defined by \eqref{eq:eq1.9} satisfying \
$|\omega |\geq 1$. \ Thus, we have
\begin{equation}
   \label{eq:eq3.2}
\lim _{n\to\infty} n^{-2}{r_n}^{-1}\big|\gamma _{r_n}^2-\delta
_{r_n}^2\big | ^{1/2}
{\mathsf E}B_n  =\Theta _{\alpha, \beta,\omega}.
\end{equation}
Observe, that \ $\lim_{n\to\infty} D_{r_n}=\lim_{n\to\infty}
\theta(\alpha,\beta,\omega_{r_n})$. 

By Remark \ref{rem:rem5} in the remaining part of the proof we may
assume \  $\alpha_{r_n}\geq0, \ \beta_{r_n}\geq0$.  \ Hence,
using Lemma \ref{lem:lem3} we have
\begin{equation}
   \label{eq:eq3.3}
\var \Bigg ( \sum _{(i,j)\in T_{k_n,\ell _n}}\big
(X_{i-1,j}^{(r_n)}\big )^2\Bigg )
\leq 2M_4\!\!\!\!\!\!\sum _{(i_1,j_1)\in T_{k_n,\ell _n}}
\sum _{(i_2,j_2)\in T_{k_n,\ell _n}}\!\!\!\!\!\!\cov \Big (
X_{i_1-1,j_1}^{(r_n)}, X_{i_2-1,j_2}^{(r_n)} \Big ) ^2\!, 
\end{equation}
where \ $M_4:=\sup_{n\in {\mathbb N}}{\mathsf
E}(\varepsilon_{0,0}^{(n)})^4$, \ and from the stationarity of \ \
$\big\{ X_{k,\ell}^{(r_n)}: k,\ell \in {\mathbb Z}\big\}$ \  follows
that the triangle \ $T_{k_n,\ell _n}$ \ can be replaced by \
$T_{n,0}$.

Now,  \eqref{eq:eq3.3} implies that if
\ $0<|\alpha| <1$ \ and \ $|\beta|=1-|\alpha|$
\begin{align}
   \label{eq:eq3.4}
n^{-4}r_n^{-1}&\big (|\gamma_{r_n}|+|\delta _{r_n}|\big) \var \Bigg (
\sum _{(i,j)\in T_{k_n,\ell _n}}\big (X_{i-1,j}^{(r_n)}\big )^2\Bigg ) \\
&\leq 2M_4 \iint\limits_T\iint\limits_T \Big( \big
(|\gamma_{r_n}|+|\delta _{r_n}|\big)^{1/2}  \cov \big
(Z_{0,1}(s_1,t_1),Z_{0,1}(s_2,t_2) \big) \Big )^2 {\mathrm
  d}s_1{\mathrm d}t_1 {\mathrm d}s_2 {\mathrm d}t_2, \nonumber
\end{align}
while for \ $|\alpha|\in\{0,1\}$,  \ $|\beta|=1-|\alpha|$  \ we have
\begin{align}
    \label{eq:eq3.5}
n^{-4}r_n^{-2}&\big |\gamma_{r_n}^2-\delta _{r_n}^2\big| \var \Bigg (
\sum _{(i,j)\in T_{k_n,\ell _n}}\big (X_{i-1,j}^{(r_n)}\big )^2\Bigg ) \\
&\leq 2M_4 \iint\limits_T\iint\limits_T \Big( \big
|\gamma_{r_n}^2-\delta _{r_n}^2\big|^{1/2}  \cov \big
(Y_{0,1}(s_1,t_1),Y_{0,1}(s_2,t_2) \big) \Big )^2 {\mathrm
  d}s_1{\mathrm d}t_1 {\mathrm d}s_2 {\mathrm d}t_2, \nonumber
\end{align}
where \ $T:=\big\{ (s,t)\in {\mathbb R}^2: 0\leq s\leq 1, -s\leq t\leq
0\big\}$. \ As the area of the triangle \ $T$ \ is finite and the
integrands in both cases are uniformly bounded on \ $T\times T$, \
Fatou's lemma and Proposition \ref{covlim} imply that the right hand
sides of \eqref{eq:eq3.4} and \eqref{eq:eq3.5} converge to \ $0$ \ as
\ $n\to \infty$. \ In a similar way one can show
\begin{equation*}
n^{-4}\kappa_n\var \Bigg ( \sum _{(i,j)\in T_{k_n,\ell_n}} \!\!\!\!\! 
X_{i-1,j}^{(r_n)}X_{i,j-1}^{(r_n)}\Bigg )\to 0 \quad \text{and} \quad
n^{-4}\kappa _n\var \Bigg ( \sum _{(i,j)\in T_{k_n,\ell_n}}\!\!\!\!\! 
\big(X_{i,j-1}^{(r_n)}\big)^2\Bigg )\to 0,
\end{equation*}
as \ $n\to \infty$, \ where
\begin{equation}
   \label{eq:eq3.6}
\kappa_n=\begin{cases}
          r_n^{-1}\big (|\gamma_{r_n}|+|\delta_{r_n}|\big) & \text{if \
            $0<|\alpha| <1$, \ $|\beta|=1-|\alpha|$},  \\
          r_n^{-2}\big |\gamma_{r_n}^2-\delta_{r_n}^2\big| & \text{if
            \ $|\alpha|\in\{0,1\}$,  \ $|\beta|=1-|\alpha|$}.
         \end{cases}
\end{equation}
that finishes the proof of Proposition \ref{pro:pro1}. \proofend

\section{Proof of Proposition \ref{pro:pro2}}
   \label{sec:sec4}
To prove Proposition \ref{pro:pro2} we are going to use the same
technique as in \cite{bpz2,bpz3}. For a given \ $n\in{\mathbb N}$ \
and \ $1\leq m\leq n$, \ let
\begin{equation*}
 A_{n,m}=\begin{pmatrix}
           A_{n,m}^{(1)} \\  A_{n,m}^{(2)}
          \end{pmatrix}:=\sum_{(k,\ell)\in T_{k_m,\ell _m}}
       \begin{pmatrix}
        X_{k-1,\ell}^{(r_n)}\varepsilon_{k,\ell}^{(r_n)} \\
        X_{k,\ell-1}^{(r_n)}\varepsilon_{k,\ell}^{(r_n)}
       \end{pmatrix},
\end{equation*}
where \ $A_{n,0}:=(0,0)^{\top}$. \ Let \ ${\mathcal F}_m^n$ \ denote
the $\sigma$-algebra generated by the random
 variables $\big\{\varepsilon^{(r_n)}_{k,\ell}:(k,\ell)\in
 U_{k_m,\ell_m}\big\}$. Obviously, \
 $A_{n,n}=A_n=\sum_{m=1}^n(A_{n,m}-A_{n,m-1})$. \
First we show that \ $\big(A_{n,m}-A_{n,m-1},{\mathcal F}_m^n\big)$ \
is a square integrable martingale difference.
Let \ $R_m:=T_{k_m,\ell_m}\setminus T_{k_{m-1},\ell_{m-1}}$, \ where \
\ $R_1:=T_{k_1,\ell_1}$. \ Short calculation shows 
\begin{equation}
    \label{eq:eq4.1}
  A_{n,m}-A_{n,m-1}
  =A_{n,m,1}+\sum_{(k,\ell)\in R_m}
            \varepsilon^{(r_n)}_{k,\ell}A_{n,m,2,k,\ell},
 \end{equation}
 where \ $A_{n,m,1}=\big(A_{n,m,1}^{(1)},A_{n,m,1}^{(2)}\big)^{\top}$ \ and \
$A_{n,m,2,k,\ell}=\big (\widetilde A_{n,m,2,k-1,\ell},\widetilde
A_{n,m,2,k,\ell-1 } \big )^{\top}$ \ with
 \begin{align}
  A_{n,m,1}^{(1)}
  &:=\sum_{(k,\ell)\in R_m}\varepsilon_{k,\ell}^{(r_n)}
                \sum_{(i,j)\in U_{k-1,\ell}\setminus U_{k_{m-1},\ell _{m-1}}}
                 \binom{k+\ell-1-i-j}{k-1-i}\alpha^{k-1-i}_{r_n}
                  \beta^{\ell-j}_{r_n}\varepsilon^{(r_n)}_{i,j}, 
\nonumber \\
 A_{n,m,1}^{(2)}
  &:=\sum_{(k,\ell)\in R_m}\varepsilon_{k,\ell}^{(r_n)}
                \sum_{(i,j)\in U_{k,\ell -1}\setminus U_{k_{m-1},\ell _{m-1}}}
                 \binom{k+\ell-1-i-j}{k-i}\alpha_{r_n}^{k-i}
                  \beta_{r_n}^{\ell-1-j}\varepsilon_{i,j}^{(r_n)},
\nonumber \\
  \widetilde A_{n,m,2,k,\ell}
  &:=\sum_{(i,j)\in U_{k,\ell}\cap U_{k_{m-1},\ell_{m-1}}}
                 \binom{k+\ell-i-j}{k-i}\alpha_{r_n}^{k-i}
                  \beta_{r_n}^{\ell-j}\varepsilon_{i,j}^{(r_n)}.
\label{eq:eq4.2}
 \end{align}
We remark that for the odd values of \ $m$ \ we have \
$R_m=\bigcup\limits_{i=-\ell _m+1}^{k_m}\big\{ (i,\ell_m)\big\}$, \ and
\begin{equation}
   \label{eq:eq4.3}
A_{n,m,1}^{(1)}= \sum_{k=-\ell
  _m+2}^{k_m}\sum_{i=-\infty}^{k-1}\alpha_{r_n} 
 ^{k-1-i}\varepsilon _{k,\ell_m}^{(r_n)}\varepsilon
 _{i,\ell_m}^{(r_n)}, \qquad  \qquad A_{n,m,1}^{(2)}=0, 
\end{equation}
while for the even values \ $R_m=\bigcup\limits_{j=-k_m+1}^{\ell_m}\big\{
(k_m,j)\big\}$, \ and 
\begin{equation}
   \label{eq:eq4.4}
A_{n,m,1}^{(2)}=\sum_{\ell =-k_m+2}^{\ell _m}\sum_{j=-\infty}^{\ell
 -1}\beta_{r_n}^{\ell-1-j} \varepsilon _{k_m,\ell}^{(r_n)}\varepsilon
_{k_m,j}^{(r_n)},  \qquad  \qquad A_{n,m,1}^{(1)}=0.
\end{equation}
 
The components of \ $A_{n,m,1}$ \ are quadratic forms of the variables \
 $\big\{\varepsilon^{(r_n)}_{i,j}:(i,j)\in R_m\big\}$, \ hence \
 $A_{n,m,1}$ \ is independent of \ ${\mathcal F}_{m-1}^n$. \ Further,
the terms $\widetilde A_{n,m,2,k,\ell}$ are linear combinations of the
variables  $\{\varepsilon^{(r_n)}_{i,j}:(i,j)\in U_{k_{m-1},\ell
  _{m-1}}\}$, \ thus they are measurable with respect to ${\mathcal
  F}_{m-1}^n$. \ Hence,
\begin{equation*}
 {\mathsf E}\big (A_{n,m}-A_{n,m-1}\mid {\mathcal F}_{m-1}^n\big )
   ={\mathsf E} A_{n,m,1}
    +\sum_{(k,\ell)\in R_m}
      A_{n,m,2,k,\ell}{\mathsf E}\big
      (\varepsilon^{(r_n)}_{p,q}\mid{\mathcal F}_{m-1}^n\big )
   =0.
\end{equation*}
By the Martingale Central Limit Theorem
 (see, e.g. \cite[Theorem 4, p. 511]{sh}), the statement in
Proposition \ref{pro:pro2} is a consequence of
 the following two propositions, where $\bone_H$ denotes the indicator function
 of the set $H$.

\begin{Pro}
    \label{pro:pro5}
If \  $0<|\alpha| <1, \ |\beta|=1-|\alpha|$ \ and \eqref{eq:eq1.7} holds then
\begin{align*}
n^{-2}r_n^{-1/2}\big(|\gamma _{r_n}|+|\delta _{r_n}|\big
)^{1/2}\sum_{m=1}^n 
   {\mathsf
     E}\big((A_{n,m}-A_{n,m-1})(A_{n,m}-&A_{n,m-1})^{\top}\big|
{\mathcal F}^n_{m-1}\big)\\ &\qmean   \big (32|\alpha||\beta |\big
)^{-1/2}\Psi_{\alpha,\beta} 
\end{align*}
as \ $n\to\infty$.

If \  $0<|\alpha|\in\{0,1\}, \ |\beta|=1-|\alpha|$ \ and
\eqref{eq:eq1.8} holds then
\begin{equation*}
n^{-2}r_n^{-1}\big |\gamma _{r_n}^2-\delta _{r_n}^2\big |^{1/2}\sum_{m=1}^n
   {\mathsf
     E}\big((A_{n,m}-A_{n,m-1})(A_{n,m}-A_{n,m-1})^{\top}\big|
{\mathcal F}^n_{m-1}\big)  \qmean  \Theta_{\alpha,\beta,\omega}
\end{equation*}
as \ $n\to\infty$.
\end{Pro}

\begin{Pro}
    \label{pro:pro6}
If \  $0<|\alpha| <1, \ |\beta|=1-|\alpha|$ \ and \eqref{eq:eq1.7} holds then
for all \ $\delta>0$
\begin{align*}
n^{-2}r_n^{-1/2}\big(|\gamma _{r_n}|+|\delta _{r_n}|\big )^{1/2}
  \sum_{m=1}^n{\mathsf E}&\big (\Vert A_{n,m}-A_{n,m-1} \Vert ^2 \\
            &\times\bone_{\left\{\Vert A_{n,m}-A_{n,m-1}\Vert
    \geq\delta nr_n^{1/4}(|\gamma _{r_n}|+|\delta _{r_n}|
    )^{-1/4}\right\}} \,\big|\,{\mathcal F}^n_{m-1}\big)
\end{align*}
converges to \ $0$ \ in probability as \ $n\to\infty$.

If \  $0<|\alpha|\in\{0,1\}, \ |\beta|=1-|\alpha|$ \ and
\eqref{eq:eq1.8} holds then for all \ $\delta>0$ 
\begin{align*}
n^{-2}r_n^{-1}\big |\gamma _{r_n}^2-\delta _{r_n}^2\big |^{1/2}\sum_{m=1}^n
   {\mathsf E}&\big (\Vert A_{n,m}-A_{n,m-1} \Vert ^2 \\
            &\times\bone_{\left\{\Vert A_{n,m}-A_{n,m-1}\Vert
    \geq\delta nr_n^{1/2}|\gamma _{r_n}^2-\delta _{r_n}^2|
    ^{-1/4}\right\}} \,\big|\,{\mathcal F}^n_{m-1}\big)
\end{align*}
converges to \ $0$ \ in probability as \ $n\to\infty$.
\end{Pro}

\noindent \textbf{Proof of Proposition \ref{pro:pro5}.} \ Let \
$U_m^n:={\mathsf E}\big
((A_{n,m}-A_{n,m-1})(A_{n,m}-A_{n,m-1})^{\top}\,\big|\, {\mathcal
  F}^n_{m-1}\big)$.  \ From the definitions
of \ $A_{n,m}$ \ and \ $B_m$ \ and from the independence of the error terms \
$\varepsilon_{k,\ell}^{(r_n)}$ \ follows that
\begin{equation*}
{\mathsf E}U^n_m={\mathsf
  E}((A_{n,m}-A_{n,m-1})(A_{n,m}-A_{n,m-1})^{\top}=
{\mathsf E}B_m-{\mathsf E}B_{m-1},
\end{equation*}
where \ $B_0$ \ is the two-by-two matrix of zeros. Thus, if \
$0<|\alpha| <1$  \ then \eqref{eq:eq3.1} implies
\begin{equation*}
\lim _{n\to\infty} n^{-2}{r_n}^{-1/2}\big(|\gamma _{r_n}|+|\delta
_{r_n}|\big ) ^{1/2}\sum_{m=1}^n
{\mathsf E}U_m^n  \to  \big (32|\alpha||\beta |\big )^{-1/2}\Psi
_{\alpha, \beta}, 
\end{equation*}
while in the case \ $|\alpha|\in\{0,1\}$  from
\eqref{eq:eq3.2} we have
\begin{equation*}
\lim _{n\to\infty} n^{-2}{r_n}^{-1}\big|\gamma _{r_n}^2-\delta
_{r_n}^2\big | ^{1/2}\sum_{m=1}^n
{\mathsf E}U_m^n  \to \Theta _{\alpha, \beta,\omega},
\end{equation*}
where $\omega$ is the limit defined by \eqref{eq:eq1.9}.

Further, from the decomposition \eqref{eq:eq4.1} follows
\begin{equation}
   \label{eq:eq4.5}
U_m^n={\mathsf E}A_{n,m,1}A_{n,m,1}^{\top }+\sum _{(k,\ell)\in
  R_m}A_{n,m,2,k,\ell }A_{n,m,2,k,\ell }^{\top }.
\end{equation}
This means that to complete the proof of the proposition we have to show
\begin{align}
\lim_{n\to\infty}n^{-4}\kappa _n\var \Big
(\sum_{m=1}^n\sum_{(k,\ell)\in
  R_m}\widetilde A_{n,m,2,k-1,\ell}^2 \Big )=0 ,
\label{eq:eq4.6} \\
\lim_{n\to\infty}n^{-4}\kappa _n\var \Big
(\sum_{m=1}^n\sum_{(k,\ell)\in R_m}\widetilde
A_{n,m,2,k-1,\ell}\widetilde A_{n,m,2,k,\ell-1}\Big )= 0 ,
\label{eq:eq4.7} \\
\lim_{n\to\infty}n^{-4}\kappa _n\var \Big
(\sum_{m=1}^n\sum_{(k,\ell)\in R_m}\widetilde A_{n,m,2,k,\ell-1}^2 \Big )=0,
\label{eq:eq4.8}
\end{align}
where \ $\kappa_n$ \ is defined by \eqref{eq:eq3.6}.

Now,  consider
 \begin{align*}
  \var\bigg(\sum_{m=1}^n\sum_{(k,\ell)\in R_m}&\widetilde
  A_{n,m,2,k,\ell}^2 \bigg)\\
 & =\sum_{m_1=1}^n\sum_{(k_1,\ell_1)\in R_{m_1}}
    \sum_{m_2=1}^n\sum_{(k_2,\ell_2)\in T_{m_2}}
     \cov\big (\widetilde A_{n,m_1,2,k_1,\ell_1}^2,
          \widetilde A_{n,m_2,2,k_2,\ell_2}^2\big ).
 \end{align*}
By Remark \ref{rem:rem5} in the remaining part of the proof we may
assume \  $\alpha_{r_n}\geq0, \ \beta_{r_n}\geq0$.  \ Hence,  as by
Lemma \ref{lem:lem3}
\begin{equation*}
 \cov\big(\widetilde A_{n,m_1,2,k_1,\ell_1}^2,
\widetilde A_{n,m_2,2,k_2,\ell_2}^2\big)
    \leq2M_4\cov\big (\widetilde A_{n,m_1,2,k_1,\ell_1},
 \widetilde A_{n,m_2,2,k_2,\ell_2}\big)^2
\end{equation*}
and representation
 \eqref{eq:eq1.3} implies
\begin{equation*}
\cov\big(\widetilde
 A_{n,m_1,2,k_1-1,\ell_1},\widetilde A_{n,m_2,2,k_2-1,\ell_2}\big)
    \leq\cov\big(X_{k_1,\ell_1}^{(r_n)},X_{k_2,\ell_2}^{(r_n)}\big),
\end{equation*}
we have
\begin{equation*}
  \var\bigg(\sum_{m=1}^n\sum_{(k,\ell)\in R_m}\widetilde
  A_{n,m,2,k,\ell}^2 \bigg) \leq 2M_4\sum_{(k_1,\ell_1)\in T_{k_n,\ell_n}}
\sum_{(k_1,\ell_1)\in
  T_{k_n,\ell_n}}\cov\big
(X_{k_1,\ell_1}^{(r_n)},X_{k_2,\ell_2}^{(r_n)}\big )^2.
\end{equation*}
Thus, using \eqref{eq:eq3.4} and \eqref{eq:eq3.5} for the cases  \
$0<|\alpha|<1$ \ and \ $|\alpha|\in\{0,1\}$, \ respectively,
\eqref{eq:eq4.6} follows from Proposition \eqref{covlim}. In a similar
way one can prove \eqref{eq:eq4.7} and \eqref{eq:eq4.8}. \proofend

\bigskip
\noindent \textbf{Proof of Proposition \ref{pro:pro6}.} \ To prove
the proposition it suffices to show
\begin{equation}
    \label{eq:eq4.9}
  n^{-4}\kappa_n\sum_{m=1}^n
         {\mathsf E}\big(\Vert A_{n,m}-A_{n,m-1}\Vert
         ^4\,\big|\,{\mathcal F}_{m-1}\big)
  \stoch0
 \end{equation}
as \ $n\to\infty$, where \ $\kappa_n$ \ is defined by
\eqref{eq:eq3.6}. By the decomposition \eqref{eq:eq4.1} 
\begin{equation*}
  \Vert A_{n,m}-A_{n,m-1}\Vert ^4
  \leq2^3\Vert A_{n,m,1}\Vert ^4
      +2^3\bigg\Vert\sum_{(k,\ell)\in R_m}
                 \varepsilon_{k,\ell}^{(r_n)}A_{n,m,2,k,\ell}\bigg\Vert^4.
 \end{equation*}
As \ $A_{n,m,1}$ \ is independent from \ ${\mathcal F}_{m-1}^n$ \ we
have \  \ ${\mathsf E}\big(\Vert A_{n,m,1}\Vert ^4\,\big|\,{\mathcal
  F}_{m-1}^n\big)={\mathsf E}\Vert A_{n,m,1}\Vert ^4$, \ while the  
 measurability of \ $A_{n,m,2,k,\ell}$ \ with respect to
 \ ${\mathcal F}_{m-1}^n$ \ implies
 \begin{equation*}
  {\mathsf E}\left(\bigg\Vert \sum_{(k,\ell)\in R_m}
                  \varepsilon_{k,\ell}^{(r_n)}A_{n,m,2,k,\ell}\bigg\Vert ^4
           \,\bigg|\,{\mathcal F}_{m-1}\right)
  \leq\big((M_4-3)^++3\big)
      \left(\sum_{(k,\ell)\in R_m}\Vert 
   A_{n,m,2,k,\ell}\Vert ^2\right)^2.
 \end{equation*}
Hence, in order to prove \eqref{eq:eq4.9}, it suffices to show
 \begin{align}
  \lim_{n\to\infty}n^{-4}\kappa _n\sum_{m=1}^n{\mathsf E} 
   \Vert A_{n,m,1}\Vert ^4&=0,\label{eq:eq4.10}\\
  \lim_{n\to\infty}n^{-4}\kappa _n \sum_{m=1}^n{\mathsf E}\left(
  \sum_{(k,\ell)\in R_m}
              \Vert A_{n,m,2,k,\ell}\Vert ^2\right)^2&=0.\label{eq:eq4.11}
 \end{align}
It is easy to see that using \eqref{eq:eq4.3} and \eqref{eq:eq4.4} we obtain
\begin{equation*} 
\Vert A_{n,m,1}\Vert ^4\!\leq\!2^3\!\left(\sum_{k=-\ell _m\!+2}^{k_m}
\sum_{i=-\infty}^{k-1}\!\!\!\alpha_{r_n}^{k-1-i}\varepsilon_{k,\ell_m}^{(r_n)}
\varepsilon_{i,\ell_m}^{(r_n)}\right)^4\!+2^3\!\left(\sum_{\ell=-k_m\!+2}^{\ell_m}
\sum_{j=-\infty}^{\ell-1}\!\!\!\beta_{r_n}^{\ell-1-j}\varepsilon_{k_m,\ell}^{(r_n)}
\varepsilon_{k_m,j}^{(r_n)}\right)^4\!\!.    
\end{equation*}
Using Lemma 12 of \cite{bpz1} a short calculation shows 
\begin{equation*}
{\mathsf E}\Vert A_{n,m,1}\Vert ^4\leq \big ( (1-\alpha
_{r_n}^2)^{-1}+(1-\beta _{r_n}^2)^{-1}\big) O(m^2), \qquad \text {as \
  $n\to\infty$,}  
\end{equation*}
and as \ $\kappa_n\big ( (1-\alpha_{r_n}^2)^{-1}+(1-\beta
_{r_n}^2)^{-1}\big )$ \ is bounded we obtain \eqref{eq:eq4.10}.

Furthermore, we have 
\begin{align*}
 {\mathsf E}\left(\sum_{(k,\ell)\in R_m}
              \Vert A_{nm,2,k,\ell}\Vert ^2\right)^2
=\sum_{(i_1,j_1)\in R_m}
    \sum_{(i_2,j_2)\in R_m} 
     {\mathsf E} &\Big ((\widetilde A_{n,m,2,i_1-1,j_1}^2+
\widetilde A_{n,m,2,i_1,j_1-1}^2) \\
&\times (\widetilde A_{n,m,2,i_2-1,j_2}^2+
 \widetilde A_{n,m,2,i_2,j_2-1}^2)\Big ).
\end{align*}
>From Lemma \ref{lem:lem3} follows 
\begin{equation*}
  {\mathsf E}\big(\widetilde A_{n,m,2,i_1,j_1}^2\widetilde
  A_{n,m,2,i_2,j_2}^2\big) 
  \leq3M_4 {\mathsf E} \widetilde A_{n,m,2,i_1,j_1}^2
{\mathsf E} \widetilde A_{n,m,2,i_2,j_2}^2,
\end{equation*}
while using \eqref{eq:eq4.2} and representation \eqref{eq:eq1.3}  one can see
 \begin{equation*}
  {\mathsf E} \widetilde A_{n,m,2,k,\ell}^2
  \leq \var X_{k,\ell}=R_{0,0}.
 \end{equation*}
Thus, 
\begin{equation*}
{\mathsf E}\left(\sum_{(k,\ell)\in R_m}
              \Vert A_{n,m,2,k,\ell}\Vert ^2\right)^2
\leq 12M_4R_{0,0}^2m^2= 12M_4\sigma _{\alpha _{r_n},\beta _{r_n}}^4m^2
\end{equation*}
that together with  \eqref{eq:eq2.6} and \eqref{eq:eq2.9} implies
\eqref{eq:eq4.11}. \proofend

\section{Proof of Proposition \ref{pro:pro3}}
   \label{sec:sec5}

In what follows we will assume  \ $0<\alpha <1$ \ and \
$\beta=1-\alpha $, \ so without loss of generality we may suppose
\ $\alpha _{r_n}, \ \beta _{r_n}, \ \gamma _{r_n}$ \ and \
$\delta _{r_n}$ \ are all positive. Consider the following
expression of \ $\det B_n$
\begin{equation*}
\det B_n=\sum _{(i_1,j_1)\in T_{k_n,\ell _n}}\sum _{(i_2,j_2)\in
T_{k_n,\ell_n}} W_{i_1,j_1,i_2,j_2}^{(n)},
\end{equation*}
where
\begin{equation*}
W_{i_1,j_1,i_2,j_2}^{(n)}:=\big (X_{i_1,j_1-1}^{(r_n)}\big )^2 \big
(X_{i_2-1,j_2}^{(r_n)}\big )^2- X_{i_1-1,j_1}^{(r_n)}
 X_{i_1,j_1-1}^{(r_n)} X_{i_2-1,j_2}^{(r_n)} X_{i_2,j_2-1}^{(r_n)}.
\end{equation*}
Using representation \eqref{eq:eq1.3} \ from Lemma \ref{lem:lem3} we
obtain.
\begin{equation}
  \label{eq:eq5.1}
{\mathsf E} W_{i_1,j_1,i_2,j_2}^{(n)}=
A^{(1,n)}_{i_1,j_1,i_2,j_2}+A^{(2,n)}_{i_1,j_1,i_2,j_2}
+A^{(3,n)}_{i_1,j_1,i_2,j_2}+ A^{(4,n)}_{i_1,j_1,i_2,j_2},
\end{equation}
where
\begin{align*}
A^{(1,n)}_{i_1,j_1,i_2,j_2}\!\!:=&\!\!\!\!\!\!\!\!\!\!\!\!\sum_{(u,v)\in
U_{(i_1-1)\land i_2,j_1\land (j_2-1)} }\!\!\!\!\!\!\!\!\!\!\!\! \big
({\mathsf E} (\varepsilon _{0,0}^{(r_n)})^4\!-\!3\big )\binom
{i_1\!+\!j_1\!-\!1\!-\!u\!-\!v}{i_1\!-\!1\!-\!u}^2
\binom{i_2\!+\!j_2\!-\!1\!-\!u\!-\!v}{i_2\!-\!u}^2 \\
&\phantom{====}\times \alpha _{r_n}^{2i_1+2i_2-2-4u}
\beta _{r_n}^{2j_1+2j_2-2-4v} \\
&-\!\!\!\!\!\!\!\!\!\!\!\!\sum_{(u,v)\in U_{i_1\land i_2-1,j_1\land
j_2-1} }\!\!\!\!\!\!\!\!\!\!\!\! \big ({\mathsf E} (\varepsilon
_{0,0}^{(r_n)})^4\!-\!3\big )\binom
{i_1\!+\!j_1\!-\!1\!-\!u\!-\!v}{i_1\!-\!1\!-\!u} \binom
{i_1\!+\!j_1\!-\!1\!-\!u\!-\!v}{i_1\!-\!u}\\
&\phantom{====}\times
\binom{i_2\!+\!j_2\!-\!1\!-\!u\!-\!v}{i_2\!-\!1\!-\!u}
\binom{i_2\!+\!j_2\!-\!1\!-\!u\!-\!v}{i_2\!-\!u}\alpha
_{r_n}^{2i_1+2i_2-2-4u}\beta _{r_n}^{2j_1+2j_2-2-4v} \!\!\!\!,\\ 
A^{(2,n)}_{i_1,j_1,i_2,j_2}\!\!:=\!&\cov \!\big
(X_{i_1-1,j_1}^{(r_n)},X_{i_2,j_2-1}^{(r_n)}\big )^2 -\cov \!\big
(X_{i_1-1,j_1}^{(r_n)},X_{i_2,j_2-1}^{(r_n)}\big )\cov \!\big
(X_{i_1,j_1-1}^{(r_n)},X_{i_2-1,j_2}^{(r_n)}\big ), \\ 
A^{(3,n)}_{i_1,j_1,i_2,j_2}\!\!:=&\!\cov \!\big
(X_{i_1-1,j_1}^{(r_n)},X_{i_2,j_2-1}^{(r_n)}\big )^2 -\cov \!\big
(X_{i_1-1,j_1}^{(r_n)},X_{i_2-1,j_2}^{(r_n)}\big )\cov \!\big
(X_{i_1,j_1-1}^{(r_n)},X_{i_2,j_2-1}^{(r_n)}\big ),\\
A^{(4,n)}_{i_1,j_1,i_2,j_2}\!\!:=&\!\var \!\big(X_{i_1-1,j_1}^{(r_n)}\big
)\!\var \!\big (X_{i_2,j_2-1}^{(r_n)}\big )\!-\!\cov\! \big
(X_{i_1-1,j_1}^{(r_n)},X_{i_1,j_1-1}^{(r_n)}\big )\! \cov\! \big
(X_{i_2-1,j_2}^{(r_n)},X_{i_2,j_2-1}^{(r_n)}\big) \\
=&\sigma _{\alpha _{r_n},\beta _{r_n}}^2\frac {1-\alpha
_{r_n}^2-\beta _{r_n}^2-\sigma _{\alpha _{r_n},\beta
_{r_n}}^{-2}}{2\alpha _{r_n}^2\beta _{r_n}^2}.
\end{align*}
Short calculation shows
\begin{equation*}
\big |A^{(1,n)}_{i_1,j_1,i_2,j_2}\big |
\leq 2(M_4+3)\cov \big
(X_{i_1-1,j_1}^{(r_n)},X_{i_2,j_2-1}^{(r_n)}\big),
\end{equation*}
so we have
\begin{align*}
n^{-4}r_n^{-1/2}&\big (\gamma _{r_n}+\delta _{r_n}\big )^{1/2}\sum
_{(i_1,j_1)\in T_{k_n,\ell _n}}\sum _{(i_2,j_2)\in T_{k_n,\ell _n}}
\big |A^{(1,n)}_{i_1,j_1,i_2,j_2}\big | \\
&\leq 2(M_4+3)\iint\limits_T\iint\limits_T \big (\gamma_{r_n}+\delta
_{r_n}\big)^{1/2}  \cov \big (Z_{0,1}(s_1,t_1),Z_{0,1}(s_2,t_2)
\big) {\mathrm d}s_1{\mathrm d}t_1 {\mathrm d}s_2 {\mathrm d}t_2.
\end{align*}
Hence, using the same arguments as in the proof of Proposition
\ref{pro:pro1} Fatou's lemma and Proposition \ref{covlim} imply
\begin{equation}
    \label{eq:eq5.2}
\lim_{n\to\infty} n^{-4}r_n^{-1/2}\big (\gamma _{r_n}+\delta
_{r_n}\big )^{1/2}\sum _{(i_1,j_1)\in T_{k_n,\ell _n}}\sum
_{(i_2,j_2)\in T_{k_n,\ell _n}}A^{(1,n)}_{i_1,j_1,i_2,j_2}=0.
\end{equation}
Next consider \
$A^{(2,n)}_{i_1,j_1,i_2,j_2}=A^{(2,n,1)}_{i_1,j_1,i_2,j_2}+
A^{(2,n,2)}_{i_1,j_1,i_2,j_2}$,  
\ where
\begin{align*}
A^{(2,n,1)}_{i_1,j_1,i_2,j_2}\!:=\!&\cov \big
(X_{i_1-1,j_1}^{(r_n)},X_{i_2,j_2-1}^{(r_n)}\big )\Big(\!\cov \big
(X_{i_1-1,j_1}^{(r_n)},X_{i_2,j_2-1}^{(r_n)}\big )\!-\!\cov \big
(X_{i_1-1,j_1}^{(r_n)},X_{i_2-1,j_2}^{(r_n)}\big )\Big ), \\
A^{(2,n,2)}_{i_1,j_1,i_2,j_2}\!:=\!&\cov \big
(X_{i_1-1,j_1}^{(r_n)},X_{i_2,j_2-1}^{(r_n)}\big )\Big(\!\cov \big
(X_{i_1-1,j_1}^{(r_n)},X_{i_2-1,j_2}^{(r_n)}\big )\!-\!\cov \big
(X_{i_1,j_1-1}^{(r_n)},X_{i_2-1,j_2}^{(r_n)}\big )\Big ).
\end{align*}
With the help of Proposition \ref{diffbound} we  can easily show
\begin{align*}
n^{-4}&r_n^{-1/2}\big (\gamma _{r_n}\!+\!\delta _{r_n}\big
)^{1/2}\sum _{(i_1,j_1)\in T_{k_n,\ell _n}}\sum _{(i_2,j_2)\in
T_{k_n,\ell _n}} \big |A^{(2,n,1)}_{i_1,j_1,i_2,j_2}\big | \\
&\leq\iint\limits_T\iint\limits_T \big (\gamma_{r_n}\!+\!\delta
_{r_n}\big)^{1/2} \Big| \cov \big (Z_{0,1}(s_1,t_1),Z_{0,1}(s_2,t_2)
\big)\Big |  \\
&\phantom{===}\times\Big| R_{[ns_1]-[ns_2]-1,[nt_1]-[nt_2]+1}-
R_{[ns_1]-[ns_2],[nt_1]-[nt_2]} \Big |{\mathrm d}s_1{\mathrm d}t_1
{\mathrm d}s_2 {\mathrm d}t_2 \\
&\leq \frac K{(\alpha _{r_n}\beta
_{r_n})^{3/2}}\iint\limits_T\!\!\iint\limits_T \!\!\big
(\gamma_{r_n}\!+\!\delta _{r_n}\big)^{1/2} \Big| \cov \big
(Z_{0,1}(s_1,t_1),Z_{0,1}(s_2,t_2) \big)\Big | {\mathrm
d}s_1{\mathrm d}t_1 {\mathrm d}s_2 {\mathrm d}t_2 \to 0
\end{align*}
as \ $n\to\infty$. \ Naturally, the same result can be proved for \
$A^{(2,n,1)}_{i_1,j_1,i_2,j_2}$, \ so we have
\begin{equation}
    \label{eq:eq5.3}
\lim_{n\to\infty} n^{-4}r_n^{-1/2}\big (\gamma _{r_n}+\delta
_{r_n}\big )^{1/2}\sum _{(i_1,j_1)\in T_{k_n,\ell _n}}\sum
_{(i_2,j_2)\in T_{k_n,\ell _n}}A^{(2,n)}_{i_1,j_1,i_2,j_2}=0.
\end{equation}
Using similar arguments one can also prove
\begin{equation}
    \label{eq:eq5.4}
\lim_{n\to\infty} n^{-4}r_n^{-1/2}\big (\gamma _{r_n}+\delta
_{r_n}\big )^{1/2}\sum _{(i_1,j_1)\in T_{k_n,\ell _n}}\sum
_{(i_2,j_2)\in T_{k_n,\ell _n}}A^{(3,n)}_{i_1,j_1,i_2,j_2}=0.
\end{equation}
Further, as \ $k_n+\ell_n=n$, \ and \ $A^{(4,n)}_{i_1,j_1,i_2,j_2}$ \
does not depend on \ $i_1,j_1,i_2,j_2$ \  using \eqref{eq:eq2.6} we obtain
\begin{equation}
    \label{eq:eq5.5}
\lim_{n\to\infty} n^{-4}r_n^{-1/2}\big (\gamma _{r_n}+\delta
_{r_n}\big )^{1/2}\sum _{(i_1,j_1)\in T_{k_n,\ell _n}}\sum
_{(i_2,j_2)\in T_{k_n,\ell _n}}A^{(4,n)}_{i_1,j_1,i_2,j_2}=\frac
2{(8\alpha\beta)^{3/2}}.
\end{equation}
Finally, the combination of representation \eqref{eq:eq5.1} and
limits \eqref{eq:eq5.2}--\eqref{eq:eq5.5} yields
\begin{equation*}
\lim_{n\to\infty} n^{-4}r_n^{-1/2}\big (\gamma _{r_n}+\delta
_{r_n}\big )^{1/2}{\mathsf E}\det B_n=\frac 2{(8\alpha\beta)^{3/2}}.
\end{equation*}

Now, let us deal with the variance of \ $\det B_n$. \ Short
calculation shows
\begin{align}
     \label{eq:eq5.6}
&n^{-8}r_n^{-1}\big (\gamma_{r_n}+\delta_{r_n}\big ) 
\var \big (\det B_n\big ) \\
&=\frac {\gamma_{r_n}+\delta_{r_n}}{n^8 r_n} \!\!\!
\sum _{(i_1,j_1)\in T_{k_n,\ell_n}}
\sum _{(i_2,j_2)\in T_{k_n,\ell_n}}\sum _{(i_3,j_3)\in T_{k_n,\ell_n}}
\sum _{(i_4,j_4)\in T_{k_n,\ell_n}}\!\!\! \cov \big (W^{(n)}_{i_1,j_1,i_2,j_2},
W^{(n)}_{i_3,j_3,i_4,j_4} \big ) \nonumber\\
&=\!\!\iint\limits _T \!\! \iint\limits _T \!\! \iint\limits _T
\!\!\iint\limits _T 
\!\!\big (\gamma_{r_n}\!+\!\delta_{r_n}\big )\Big
(\!\Theta_n^{(1)}(s_1,t_1,s_2,t_2,s_3,t_3,s_4,t_4)\!+\! 
\Theta_n^{(2)}(s_1,t_1,s_2,t_2,s_3,t_3,s_4,t_4)  \nonumber \\ 
&\phantom{==============}+
2\Theta_n^{(3)}(s_1,t_1,s_2,t_2,s_3,t_3,s_4,t_4)\Big )
{\mathrm d}s_1{\mathrm d}t_1
{\mathrm d}s_2{\mathrm d}t_2{\mathrm d}s_3{\mathrm d}t_3{\mathrm d}s_4
{\mathrm d}t_4, \nonumber
\end{align}
where
\begin{align*}
&\Theta_n^{(1)}(s_1,t_1,s_2,t_2,s_3,t_3,s_4,t_4) \\ &:=\! \cov \!\Big (\!
Z_{1,0}^{(r_n)}(s_1,t_1)Z_{0,1}^{(r_n)}(s_2,t_2)\!
\big(X^{(n)}_{[ns_1],[nt_1]-1}\!\!-\!X^{(r_n)}_{[ns_1]-1,[nt_1]}\big )\!
\big(X^{(n)}_{[ns_2]-1,[nt_2]}\!\!-\!X^{(r_n)}_{[ns_2],[nt_2]-1}\big ), \\
&\hskip .9 truecm Z_{1,0}^{(n)}(s_3,t_3)Z_{0,1}^{(n)}(s_4,t_4)
\big (X^{(n)}_{[ns_3],[nt_3]-1}\!\!-\!X^{(r_n)}_{[ns_3]-1,[nt_3]}\big )
\big (X^{(n)}_{[ns_4]-1,[nt_4]}\!\!-\!X^{(r_n)}_{[ns_4],[nt_4]-1}\big ) \Big ),\\
&\Theta_n^{(2)}\!(s_1,t_1,s_2,t_2,s_3,t_3,s_4,t_4)\!:=\! \cov\! \Big (\!
Z_{1,0}^{(n)}(s_1,t_1)Z_{0,1}^{(n)}(s_1,t_1)
\big (X^{(n)}_{[ns_2],[nt_2]-1}\!-\!X^{(r_n)}_{[ns_2]-1,[nt_2]}\big )^2\!\!, \\
&\hskip 3.2 truecm Z_{1,0}^{(n)}(s_3,t_3)Z_{0,1}^{(n)}(s_3,t_3)
\big (X^{(n)}_{[ns_4],[nt_4]-1}\!\!-\!X^{(r_n)}_{[ns_4]-1,[nt_4]}\big
)^2\Big ), \\ 
&\Theta_n^{(3)}\!(s_1,t_1,s_2,t_2,s_3,t_3,s_4,t_4)\!:=\! \cov\! \Big (\!
Z_{1,0}^{(n)}(s_1,t_1)Z_{0,1}^{(n)}(s_1,t_1)
\big (X^{(n)}_{[ns_2],[nt_2]-1}\!-\!X^{(r_n)}_{[ns_2]-1,[nt_2]}\big )^2\!\!, \\
&\hskip .9 truecm Z_{1,0}^{(n)}(s_3,t_3)Z_{0,1}^{(n)}(s_4,t_4)
\big (X^{(n)}_{[ns_3],[nt_3]-1}\!\!-\!X^{(r_n)}_{[ns_3]-1,[nt_3]}\big )
\big (X^{(n)}_{[ns_4]-1,[nt_4]}\!-\!X^{(r_n)}_{[ns_4],[nt_4]-1}\big ) \Big ).
\end{align*}
By representation \eqref{eq:eq1.3} the components \ $\Theta_n^{(q)},
\ q=1,2,3$, \ of the integrand in the right hand side of
\eqref{eq:eq5.6} are linear combinations of 
covariances of form 
\begin{equation}
  \label{eq:eq5.7}
  \cov (\varepsilon _{i_1,j_1} \varepsilon _{i_2,j_2}\varepsilon _{i_3,j_3}
\varepsilon _{i_4,j_4},\varepsilon _{i_5,j_5} \varepsilon _{i_6,j_6}
\varepsilon _{i_7,j_7}\varepsilon _{i_8,j_8}),
\end{equation}
where the indices \ $(i_r,j_r)\in {\mathbb Z}^2, \ r=1,2,\dots ,8$, \ run 
either on quarter planes \ $U_{[ns_q],[nt_q]-1}$ \ or on \
$U_{[ns_q]-1,[nt_q]}$, \  $q=\big [ (r+1)/2 \big ]$. \ Using the
definitions of Lemma \ref{lem:lem2} we can express the coefficients of
the linear combinations as products of \ $1/r_n$ \ and two terms of form 
\ $(\alpha_{r_n}+\beta_{r_n})^{[ns_{m_r}]+[nt_{m_r}]-1-i_r-j_r}{\mathsf
  P} \big ( S^{(\nu_{r_n})}_{[ns_{m_r}]+[nt_{m_r}]-1-i_r-j_r}= 
[ns_{m_r}]-1-i_r \big )$, \ two terms of form  \
$(\alpha_{r_n}+\beta_{r_n})^{[ns_{m_r}]+[nt_{m_r}]-1-i_r-j_r}{\mathsf
  P} \big
(S^{(\nu_{r_n})}_{[ns_{m_r}]+[nt_{m_r}]-1-i_r-j_r}=[ns_{m_r}]-i_r \big
)$ \ and  four terms of form  \
$(\alpha_{r_n}+\beta_{r_n})^{[ns_{m_r}]+[nt_{m_r}]-1-i_r-j_r}\widehat\Delta
^{(n)}_{i_r,j_r}(s_{m_r},t_{m_r})$ \ with   
\begin{align*}
\widehat\Delta ^{(n)}_{i_r,j_r}&(s_{m_r},t_{m_r}) \\
:=&{\mathsf P} \big ( 
S^{(\nu_{r_n})}_{[ns_{m_r}]+[nt_{m_r}]-1-i_r-j_r}\!=\![ns_{m_r}]\!-\!i_r \big )\!-\!
{\mathsf P} \big ( 
S^{(\nu_{r_n})}_{[ns_{m_r}]+[nt_{m_r}]-1-i_r-j_r}\!=\![ns_{m_r}]\!-\!1\!-\!i_r \big
),
\end{align*}
where \ $\nu_{r_n}=\alpha_{r_n}/(\alpha_{r_n}+\beta_{r_n})$ \ and \
$\cup_{r=1}^8 \{ m_r \}=\{1,2,3,4\}$. \ Corollary 
\ref{cor:cor1} implies that there exists a positive constant \ $C$ \
such that
\begin{equation}
  \label{eq:eq5.8}
\Big\vert \widehat\Delta ^{(n)}_{i_r,j_r}(s_{m_r},t_{m_r}) \Big\vert 
\leq \frac C{\alpha_{r_n}\beta_{r_n}([ns_{m_r}]+[nt_{m_r}]-1-i_r-j_r)}.
\end{equation}
Covariances of form \eqref{eq:eq5.7} are equal to
zero if the index sets \ $\big \{ (i_r,j_r) :  r=1,2,3,4 \big\}$ \
and \ $\big \{ (i_r,j_r)   : r=5,6,7,8 \big\}$ \ are disjoint. 
Besides the nonempty intersection of these sets, to obtain nonzero covariances
in \eqref{eq:eq5.7} for  each \ $u\in \{1,2, \dots ,8 \}$ \ there 
should exist at least one \ $v\in \{1,2, \dots ,8 \}$ \ such that \
$u\ne v$\ and \ $(i_u,j_u)=(i_v,j_v)$. \ Consider first the  case, 
when \ $\{1,2,\dots,8\}$ \  is divided into two disjoint subsets
\ $\{u_1,u_2,u_3,u_4\}$ \ and \ $\{v_1,v_2,v_3,v_4\}$, \
$(i_{u_r},j_{u_r})=(i_{v_r},j_{v_r}), \ r=1,2,3,4,$ \ holds and no
other index pairs are equal. This  configuration yields the highest
amount of terms when we express the covariances of \ $\Theta_n^{(q)},
\ q=1,2,3$. \ Expression  \eqref{eq:eq5.6}
shows that the sum of the corresponding terms of \ 
$n^{-8}r_n^{-1}\big (\gamma_{r_n}+\delta_{r_n}\big ) \var \big (\det
B_n\big )$ \ can be rewritten as the sum of terms of  form 
\begin{align*}
\iint\limits _T \! \iint\limits _T &\! \iint\limits _T \!\iint\limits _T
\Big( R_{[ns_{m_1}]-[ns_{m_2}],[nt_{m_1}]-[nt_{m_2}]}-
 R_{[ns_{m_1}]-[ns_{m_2}]\pm 1,[nt_{m_1}]-[nt_{m_2}]\mp 1} \Big )\\
\times&\Big( R_{[ns_{m_3}]-[ns_{m_4}],[nt_{m_3}]-[nt_{m_4}]}-
 R_{[ns_{m_3}]-[ns_{m_4}]\pm 1,[nt_{m_3}]-[nt_{m_4}]\mp 1} \Big )\\
\times&\big (\gamma_{r_n}+\delta_{r_n}\big )^{1/2}\cov \big
(Z_{u_1,v_1}^{(n)}(s_{m_5},t_{m_5}) Z_{u_2,v_2}^{(n)}(s_{m_6},
t_{m_6})\big ) \\
\times&\big(\gamma_{r_n}+\delta_{r_n}\big )^{1/2}\cov \big
(Z_{u_3,v_3}^{(n)}(s_{m_7},t_{m_3}) Z_{u_4,v_4}^{(n)}(s_{m_8},t_{m_8})\big )
{\mathrm d}s_1{\mathrm d}t_1
{\mathrm d}s_2{\mathrm d}t_2{\mathrm d}s_3{\mathrm d}t_3{\mathrm d}s_4
{\mathrm d}t_4, 
\end{align*}
where \ $\{m_r: r=1,2,\dots ,8\}=\{1,2,3,4\}$, \ $(u_r,v_r)\in 
\{(0,1),(1,0)\}, \ r=1,2,3,4$. \ 
Fatou's lemma, Lemma \ref{lem:lem1}  and Propositions \ref{diffbound}
and \ref{covlim} imply that these terms of the sum \ 
$n^{-8}r_n^{-1}\big (\gamma_{r_n}+\delta_{r_n}\big ) \var \big (\det 
B_n\big )$ \ converge to \ $0$ \ as \ $n\to\infty$.

The next case is when  \ $\{1,2,\dots,8\}$ \  is divided into three 
disjoint subsets
\ $\{u_1,u_2,u_3\}$ \ and \ $\{v_1,v_2,v_3\}$ \ and \  $\{w_1,w_2\}$ \
and either  
\begin{equation}
   \label{eq:eq5.9}
(i_{u_r},j_{u_r})=(i_{v_r},j_{v_r})=(i_{w_r},j_{w_r}), \  r=1,2, 
\qquad \text{and} \qquad (i_{u_3},j_{u_3})=(i_{v_3},j_{v_3})
\end{equation}
or
\begin{equation}
    \label{eq:eq5.10}
(i_{u_r},j_{u_r})=(i_{v_r},j_{v_r}), \ r=1,2, 
\quad \text{and} \quad
(i_{u_3},j_{u_3})=(i_{v_3},j_{v_3})=(i_{w_1},j_{w_1})=(i_{w_1},j_{w_2})
\end{equation}
holds and no other index pairs are equal. 
Inequality
\eqref{eq:eq5.8} implies  that we have
\begin{align*}
\!\!\sum _{(i,j)\in U_{[ns_1]\land [ns_2]\land [ns_3]-1, [nt_1]\land
  [nt_2]\land [nt_3]-1}}\!\!\!\!\!\!\!\!\!\!\!\!\!\!\!\!\!\!\!\! 
&{\mathsf P} \big (
S^{\alpha}_{[ns_1]+[nt_1]-1-i-j}\!=\! 
[ns_1]\!-\!i \big ) \big |\widehat\Delta ^{(n)}_{i,j}(s_2,t_2)\Big |
 \big |\widehat\Delta ^{(n)}_{i,j}(s_3,t_3)\Big | \\ 
\leq \!\!\!\!\!\sum_{m=-\infty}^{[ns_1]\land [ns_2]\land [ns_3]+[nt_1]\land
  [nt_2]\land [nt_3]-2}\!\!\!\!\!\!\!\!\!\!\!\!\!\!\!  &\frac
{C^2}{(\alpha_{r_n}\beta_{r_n})^2([ns_2]+[nt_2]-1-m)([ns_3]+[nt_3]-1-m)}
\\ 
\times\!\!\!\!\!\!\!\!\!\sum_{i=m-[nt_1]\land [nt_2]\land
  [nt_3]+1}^{[ns_1]\land [ns_2]\land [ns_3]-1}&\!\!\!\!\!\!\!\!\!
 {\mathsf P} \big ( S^{\alpha}_{[ns_1]+[nt_1]-1-m}\!=\![ns_1]\!-\!i
 \big ) \!\leq\! 
 \frac {C^2\zeta(2)}{(\alpha_{r_n}\beta_{r_n})^2}, 
\end{align*}
so the expressions of the above form are bounded uniformly in \ $n$ \
and  \ $(s_r,t_r)\in T, \ r=1,2,3$. \ 
Similarly, by Remark \ref{rem:rem4}  there exists a constant
$D>0$ such that 
\begin{equation}
  \label{eq:eq5.11}
\sum _{(i,j)\in \bigcap _{r=1}^4U_{[ns_r]+u_r,[nt_r]+v_r}}\!\!\!\! \Big (
\prod _{r=1}^4{\mathsf P} \big (S^{\alpha}_{[ns_r]+[nt_r]-1-i-j}\!=\!
[ns_r]\!-\!i\!-\!u_r \big )  \Big)  \!\leq\! \frac
{D^3\zeta(3/2)}{(\alpha_{r_n}\beta_{r_n})^{3/2}},  
\end{equation}
where \ $(u_r,v_r)\in \{(0,1),(1,0)\}, \ r=1,2,3,4$. \ 

It is not difficult to show that in the case described by \eqref{eq:eq5.9} 
the corresponding part of the sum  \ $n^{-8}r_n^{-1}\big
(\gamma_{r_n}+\delta_{r_n}\big ) \var \big (\det B_n\big )$ \ 
can always be bounded from above by the sum of components of the form
\begin{align*}
\iint\limits _T \! \iint\limits _T \! &\iint\limits _T \!\iint\limits _T
 \frac {C^2\zeta(2)}{(\alpha_{r_n}\beta_{r_n})^2} 
\big (\gamma_{r_n}+\delta_{r_n}\big )^{1/2}\Big |\cov \big
(Z_{u_1,v_1}^{(n)}(s_{m_1},t_{m_1}) Z_{u_2,v_2}^{(n)}(s_{m_2},
t_{m_2})\big )\Big| \\
\times&\big(\gamma_{r_n}+\delta_{r_n}\big )^{1/2}\Big |\cov \big
(Z_{u_3,v_3}^{(n)}(s_{m_3},t_{m_3})
Z_{u_4,v_4}^{(n)}(s_{m_4},t_{m_4})\big )\Big |
{\mathrm d}s_1{\mathrm d}t_1
{\mathrm d}s_2{\mathrm d}t_2{\mathrm d}s_3{\mathrm d}t_3{\mathrm d}s_4
{\mathrm d}t_4 
\end{align*}
where \ $\{m_r: r=1,2,3,4\}\subseteq \{1,2,3,4\}$ contains at least
$3$ different points and \ $(u_r,v_r)\in 
\{(0,1),(1,0)\}, \ r=1,2,3,4$. \ In this way  by Fatou's lemma
and Proposition \ref{covlim} we obtain that the terms of
 \ $n^{-8}r_n^{-1}\big
(\gamma_{r_n}+\delta_{r_n}\big ) \var \big (\det B_n\big )$ \
corresponding to case \eqref{eq:eq5.9} converge
to \ $0$ \ 
as \ $n\to\infty$. \ Using similar ideas and \eqref{eq:eq5.11} the
same can be proved in the case \eqref{eq:eq5.10}. 

The remaining terms of  \ $n^{-8}r_n^{-1}\big
(\gamma_{r_n}+\delta_{r_n}\big ) \var \big (\det B_n\big )$ \ can be
handled in a similar way. \proofend

\section{Proof of Proposition \ref{pro:pro4}}
   \label{sec:sec6}

Similarly to Section \ref{sec:sec5} it is enough to consider the case 
\ $0<\alpha <1$ \ and \ $\beta=1-\alpha $. \ 
We have
\begin{align*}
n^{-3}r_n^{-1/2}\big (\gamma _{r_n}+\delta
_{r_n}\big )^{1/2}\bar B_n A_n=\Big (n^{-2}r_n^{-1/2}\big (\gamma _{r_n}+\delta
_{r_n}\big )^{1/2}\bar B_n&-\frac 1{\sqrt{32\alpha\beta }}\bar{\mathsf 1}\Big
)\frac 1nA_n \\&+\frac 1{\sqrt{32\alpha\beta }}\frac 1n\bar{\mathsf 1}A_n, 
\end{align*}
where \ ${\mathsf 1}$ \ denotes the two-by-two matrix of ones. Short
straightforward calculations shows 
\begin{equation*}
\Big (n^{-2}r_n^{-1/2}\big (\gamma _{r_n}+\delta
_{r_n}\big )^{1/2}\bar B_n-\frac 1{\sqrt{32\alpha\beta }}\bar{\mathsf 1}\Big
)\frac 1n A_n=C_n+D_n,
\end{equation*}
where
\begin{align*}
C_n&:=n^{-1}r_n^{-1/4}\big (\gamma _{r_n}+\delta_{r_n}\big
)^{1/4}\diag(A_n)n^{-2}r_n^{-1/4}\big (\gamma _{r_n}+\delta_{r_n}\big
)^{1/4} \bar B_n(1,1)^{\top},\\
D_n&:=\Big (n^{-2}r_n^{-1/2}\big (\gamma _{r_n}+\delta
_{r_n}\big )^{1/2}\sum_{(i,j)\in
  T_{k_n,\ell_n}}X_{i-1,j}^{(r_n)}X_{i,j-1}^{(r_n)} - 
\frac 1{\sqrt{32\alpha\beta }}\Big ) 
\frac 1n Q_n(1,-1)^{\top}.
\end{align*}
Here \ $\diag (A_n)$ \ denotes the two-by-two diagonal matrix
having \ $A_n$ \ in its main diagonal and
\begin{equation}
    \label{eq:eq6.1}
Q_n:=(1,-1)A_n=\sum_{(i,j)\in  T_{k_n,\ell_n}}
\big (X_{i-1,j}^{(r_n)}-X_{i,j -1}^{(r_n)}\big )\varepsilon
_{i,j }. 
\end{equation} 
By Proposition \ref{pro:pro1}
\begin{equation}
    \label{eq:eq6.2}
n^{-2}r_n^{-1/2}\big (\gamma _{r_n}+\delta
_{r_n}\big )^{1/2}\sum_{(i,j)\in
  T_{k_n,\ell_n}}X_{i-1,j}^{(r_n)}X_{i,j-1}^{(r_n)} - 
\frac 1{\sqrt{32\alpha\beta }}
\qmean 0 \qquad \text{as \ 
$n\to\infty$.}
\end{equation}
Representation \eqref{eq:eq1.3} and independence of
the error terms \ $\varepsilon _{i,j}^{(r_n)}$ \ imply \ ${\mathsf E}
Q_n=0$ \ and 
\begin{align*}
{\mathsf E}Q_n^2=&\sum_{(i,j)\in  T_{k_n,\ell_n}}
{\mathsf E}\big (X_{i-1,j}^{(r_n)}-X_{i,j -1}^{(r_n)}\big
)^2=(k_n+\ell_n)(k_n+\ell _n+1)\big (R_{0,0}-R_{-1,1}\big ) \\
=&\frac {n(n+1)}{4\alpha _{r_n}\beta _{r_n}}\bigg (1+
\Big (\frac {\gamma _{r_n}+\delta_{r_n}}{r_n}\Big)^{1/2}
\sigma _{\alpha _{r_n},\beta_{r_n}}^2
\Big (\frac {\gamma _{r_n}+\delta_{r_n}}{r_n}\Big)^{1/2}
\Big (\frac {\gamma _{r_n}+\delta_{r_n}}{r_n}-2\Big)
\bigg ).
\end{align*}
Taking into account \eqref{eq:eq2.6} we obtain \ 
\begin{equation}
   \label{eq:eq6.3}
\lim_{n\to\infty}\frac 1{n^2}{\mathsf E}
Q_n^2= \frac 1{4\alpha\beta}
\end{equation}
that together with \eqref{eq:eq6.2}
implies \ $D_n\stoch (0,0)^{\top}$ \ as \  $n\to\infty$. 
\begin{align*}
\frac 1{n^2}\Big (&\frac {\gamma _{r_n}\!+\!\delta _{r_n}}{r_n}\Big
)^{1/4}\!{\mathsf E}\Big (\bar B_n (1,1)^{\top} \Big )\!=\!\frac 1{n^2}
\Big (\frac {\gamma _{r_n}\!+\!\delta _{r_n}}{r_n}\Big )^{1/4} 
\!\!\!\!\!\!\!\sum_{(k,\ell)\in T_{k_n,\ell_n}} \!\!\!\!\!\!{\mathsf E}
       \begin{pmatrix}
        \big (X_{k,\ell-1}^{(r_n)}\big)^2\!-\!X_{k-1,\ell}^{(r_n)}X_{k,\ell
        -1}^{(r_n)} \\ 
        \big(X_{k,\ell-1}^{(r_n)}\big)^2\!-\!X_{k-1,\ell}^{(r_n)}X_{k,\ell
        -1}^{(r_n)}  
       \end{pmatrix} \\
&= \Big (\frac {\gamma _{r_n}\!+\!\delta _{r_n}}{r_n}\Big
)^{1/4}\frac{n+1}{2n}\big (R_{0,0}\!-\!R_{-1,1}\big )\begin{pmatrix}1\\1
       \end{pmatrix}\!=\!
       \Big (\frac {\gamma _{r_n}\!+\!\delta _{r_n}}{r_n}\Big 
)^{1/4}\frac 1{2n^2} {\mathsf E}Q_n^2\begin{pmatrix}1\\1
       \end{pmatrix}\to \begin{pmatrix}0\\0
       \end{pmatrix}
\end{align*}
as \ $n\to\infty$. \ Furthermore, with the help of Lemma \ref{lem:lem3} 
we obtain
\begin{align*}
\var \bigg (&\sum_{(k,\ell)\in T_{k_n,\ell_n}} 
\Big( \big (X_{k,\ell-1}^{(r_n)}\big)^2-X_{k-1,\ell}^{(r_n)}X_{k,\ell
  -1}^{(r_n)}\Big)\bigg )\\ &= \sum_{(i_1,j_1)\in T_{k_n,\ell_n}} 
\sum_{(i_2,j_2)\in T_{k_n,\ell_n}}
B^{(1,n)}_{i_1,j_1,i_2,j_2} +2B^{(2,n)}_{i_1,j_1,i_2,j_2} 
+B^{(3,n)}_{i_1,j_1,i_2,j_2}+B^{(4,n)}_{i_1,j_1,i_2,j_2}, 
\end{align*}
where
\begin{align*}
B^{(1,n)}_{i_1,j_1,i_2,j_2}\!:=\!&\!\!\!\!\!\!\!\!\!\!\!\!\sum_{(u,v)\in
U_{i_1\land i_2,j_1\land j_2-1} }\!\!\!\!\!\!\!\!\!\!\!\! \big
({\mathsf E} (\varepsilon _{0,0}^{(r_n)})^4\!-\!3\big )\binom
{i_1\!+\!j_1\!-\!1\!-\!u\!-\!v}{i_1\!-\!u}^2
\binom{i_2\!+\!j_2\!-\!1\!-\!u\!-\!v}{i_2\!-\!u}^2 \\
&\phantom{======}\times \alpha _{r_n}^{2i_1+2i_2-4u}
\beta _{r_n}^{2j_1+2j_2-4-4v} \\
&-\!\!\!\!\!\!\!\!\!\!\!\!\sum_{(u,v)\in U_{i_1\land (i_2-1),j_1\land
j_2-1} }\!\!\!\!\!\!\!\!\!\!\!\! 2\big ({\mathsf E} (\varepsilon
_{0,0}^{(r_n)})^4\!-\!3\big )\binom
{i_1\!+\!j_1\!-\!1\!-\!u\!-\!v}{i_1\!-\!u}^2 \binom
{i_2\!+\!j_2\!-\!1\!-\!u\!-\!v}{i_2\!-\!1\!-\!u}\\
&\phantom{====}\times\!
\binom{i_2\!+\!j_2\!-\!1\!-\!u\!-\!v}{i_2\!-\!u}\alpha
_{r_n}^{2i_1+2i_2-1-4u}\beta _{r_n}^{2j_1+2j_2-3-4v}\\
&+\!\!\!\!\!\!\!\!\!\!\!\!\sum_{(u,v)\in U_{i_1\land i_2-1,j_1\land
j_2-1} }\!\!\!\!\!\!\!\!\!\!\!\! \big ({\mathsf E} (\varepsilon
_{0,0}^{(r_n)})^4\!-\!3\big )\binom
{i_1\!+\!j_1\!-\!1\!-\!u\!-\!v}{i_1\!-\!1\!-\!u} \binom
{i_1\!+\!j_1\!-\!1\!-\!u\!-\!v}{i_1\!-\!u}\\
&\phantom{====}\times\!
\binom{i_2\!+\!j_2\!-\!1\!-\!u\!-\!v}{i_2\!-\!1\!-\!u}
\binom{i_2\!+\!j_2\!-\!1\!-\!u\!-\!v}{i_2\!-\!u}\alpha
_{r_n}^{2i_1+2i_2-2-4u}\beta _{r_n}^{2j_1+2j_2-2-4v}, \\ 
B^{(2,n)}_{i_1,j_1,i_2,j_2}\!:=\!&\cov \big
(X_{i_1,j_1-1}^{(r_n)},X_{i_2,j_2-1}^{(r_n)}\big )
\Big (\!\cov \big (X_{i_1,j_1-1}^{(r_n)},X_{i_2,j_2-1}^{(r_n)}\big
)\!-\!\cov \big(X_{i_1,j_1-1}^{(r_n)}, X_{i_2-1,j_2}^{(r_n)}\big )\Big ),\\
B^{(3,n)}_{i_1,j_1,i_2,j_2}\!:=\!&\cov \big
(X_{i_1,j_1-1}^{(r_n)},X_{i_2,j_2-1}^{(r_n)}\big )
\Big (\!\cov \big (X_{i_1-1,j_1}^{(r_n)},X_{i_2-1,j_2}^{(r_n)}\big
)\!-\!\cov (X_{i_1,j_1-1}^{(r_n)}, X_{i_2-1,j_2}^{(r_n)}\big )\Big ),\\
B^{(4,n)}_{i_1,j_1,i_2,j_2}\!:=\!&\cov \big
(X_{i_1,j_1-1}^{(r_n)},X_{i_2-1,j_2}^{(r_n)}\big ) 
\Big (\!\cov \big(X_{i_1-1,j_1}^{(r_n)},X_{i_2,j_2-1}^{(r_n)}\big )\!-\!\cov
\big(X_{i_1,j_1-1}^{(r_n)}, X_{i_2,j_2-1}^{(r_n)}\big )\Big ).
\end{align*}
Hence, using the same arguments as in the proof of Proposition
\ref{pro:pro3} (see \eqref{eq:eq5.2} and \eqref{eq:eq5.3}) one can verify 
\begin{equation*}
\lim_{n\to\infty}\frac 1{n^4}\Big (\frac {\gamma _{r_n}+\delta
  _{r_n}}{r_n}\Big )^{1/2}\var \bigg (\sum_{(k,\ell)\in T_{k_n,\ell_n}} 
\Big( \big (X_{k,\ell-1}^{(r_n)}\big)^2-X_{k-1,\ell}^{(r_n)}X_{k,\ell
  -1}^{(r_n)}\Big)\bigg )=0.
\end{equation*}
Naturally, the same holds for the second component of
\ $n^{-2}r_n^{-1/4}\big (\gamma_{r_n}+\delta _{r_n}\big )^{1/4}\bar
B_n (1,1)^{\top}$, \ that means  
\begin{equation}
   \label{eq:eq6.4}
n^{-2}r_n^{-1/4}\big (\gamma_{r_n}+\delta _{r_n}\big )^{1/4}\bar
B_n (1,1)^{\top}\qmean (0,0)^{\top} \qquad \qquad \text{as \ 
$n\to\infty$.}
\end{equation}
Proposition \ref{pro:pro2} and \eqref{eq:eq6.4} imply  \ 
$C_n\stoch (0,0)^{\top}$ \ as \ $n\to\infty$, \ so to prove the asymptotic
normality of \ $n^{-3}n^{-1/2}\big (\gamma_{r_n}+\delta _{r_n}\big
)^{1/2}\bar B_nA_n$ \ it suffices to show the asymptotic normality of
\ $n^{-1}\bar{\mathbf 1}A_n=n^{-1}Q_n (1, -1)^{\top}$.

For a given \ $n\in{\mathbb N}$ \ and \ $1\leq m\leq n$ \ let \ 
$Q_{n,m}:=(1,-1)A_{n,m}$. \ Obviously \ $Q_{n,n}=Q_n$ \ and  from
\eqref{eq:eq4.1} we have 
\begin{equation}
    \label{eq:eq6.5}
  Q_{n,m}-Q_{n,m-1}
  =A_{n,m,1}^{(1)}-A_{n,m,1}^{(2)}+\sum_{(k,\ell)\in R_m}
            \varepsilon_{k,\ell}^{(r_n)}\big (\widetilde A_{n,m,2,k-1,\ell}-
\widetilde A_{n,m,2,k,\ell -1}\big ).
 \end{equation} 
As  \ $\big(Q_{n,m}-Q_{n,m-1},{\mathcal F}_m^n\big)$ \
is a square integrable martingale difference, similarly to the proof
of Proposition \ref{pro:pro2} the statement of Proposition
\ref{pro:pro4} follows from the propositions below.

\begin{Pro}
    \label{pro:pro7}
If \  $0<\alpha <1, \ \beta=1-\alpha$ \ and \eqref{eq:eq1.7} holds then
\begin{equation*}
\frac 1{n^2}\sum_{m=1}^n
   {\mathsf
     E}\big((Q_{n,m}-Q_{n,m-1})^2\big|
{\mathcal F}^n_{m-1}\big)  \stoch  \frac 1{4\alpha\beta} \qquad \qquad
\text{as \ $n\to\infty$.}
\end{equation*}
\end{Pro}

\begin{Pro}
    \label{pro:pro8}
If \  $0<\alpha <1, \ \beta=1-\alpha$ \ and \eqref{eq:eq1.7} holds then
for all \ $\delta>0$
\begin{equation*}
\frac 1{n^2}\sum_{m=1}^n{\mathsf E}\big (( Q_{n,m}-Q_{n,m-1} )^2
            \bone_{\left\{|Q_{n,m}-Q_{n,m-1}|
    \geq\delta n\right\}} \,\big|\,{\mathcal F}^n_{m-1}\big)\stoch 0
\qquad \qquad \text{as \ $n\to\infty$.}
\end{equation*}
\end{Pro}

\medskip
\noindent
\textbf{Proof of Proposition \ref{pro:pro7}.} \ The proof is very
similar to that of Proposition \ref{pro:pro5}. Let 
$V_m^n:={\mathsf E}\big((Q_{n,m}-Q_{n,m-1})^2\,\big|\,
{\mathcal F}^n_{m-1} \big )$. \ The statement of Proposition
\ref{pro:pro7} will follow from 
 \begin{equation}
   \label{eq:eq6.6}
 \lim_{n\to\infty} \frac 1{n^2}\sum_{m=1}^n{\mathsf E} V_m^n=\frac
 1{4\alpha\beta}  \qquad \text{and} \qquad
  \lim_{n\to\infty}\frac 1{n^4}\var\bigg(\sum_{m=1}^nV_m^n\bigg)=0.
 \end{equation}
By the martingale property of \ $Q_{n,m}$ \ we have
\begin{align*}
\sum_{m=1}^n{\mathsf E} V_m^n=\sum_{m=1}^n\big({\mathsf
  E}Q_{n,m}^2-{\mathsf E}Q_{n,m-1}^2\big)={\mathsf E}Q_{n}^2 
\end{align*}
that together with \eqref{eq:eq6.3} implies the convergence of the
means in \eqref{eq:eq6.6}.
Furthermore, representations \eqref{eq:eq6.1} of \ $Q_n^m$ \  and
\eqref{eq:eq4.5} of  \ $U_n^m$ \ imply
\begin{equation*}
V^n_m=(1,-1)U_nU_m^n (1,-1)^{\top}={\mathsf E} \big
(A_{n,m,1}^{(1)}-A_{n,m,1}^{(2)}  
\big )^2+\!\!\!\sum_{(k,\ell)\in R_m}\!\!\!
            \big (\widetilde A_{n,m,2,k-1,\ell}-
\widetilde A_{n,m,2,k,\ell -1}\big )^2.
\end{equation*}
Using representation \eqref{eq:eq1.3}, definition \eqref{eq:eq4.2}
and Lemma \ref{lem:lem3} one can verify
\begin{align*}
\var\left(\sum_{m=1}^n V_m^n\right) &=
\var\bigg (\sum_{m=1}^n \sum_{(k,\ell)\in 
R_m}\big (\widetilde A_{n,m,2,k-1,\ell}-
\widetilde A_{n,m,2,k,\ell -1}\big )^2\bigg) \\ 
&\leq \sum_{(i_1,j_1)\in T_{k_n,\ell_n}}\sum_{(i_2,j_2)\in T_{k_n,\ell_n}}
G_{n,i_1,j_1,i_2,j_2}+H_n,
\end{align*} 
where 
\begin{equation*}
G_{n,i_1,j_1,i_2,j_2}:=\cov \Big( \big
(X_{i_1-1,j_1}^{(r_n)}-X_{i_1,j_1-1}^{(r_n)} 
\big )^2,\big (X_{i_2-1,j_2}^{(r_n)}-X_{i_2,j_2-1}^{(r_n)}
\big)^2\Big)
\end{equation*}
and \ $n^{-4}H_n\to 0$ \ as \ $n\to\infty$. \ As  
 \ $X_{k-1,\ell}^{(r_n)}-X_{k,\ell -1}^{(r_n)}$ \ is also a linear combination
of the variables \ $\big\{ \varepsilon_{i,j}^{(r_n)} : (i,j)\in
U_{k,\ell }\big \}$, \ by Lemma \ref{lem:lem3} we have
\begin{align*} 
&\sum_{(i_1,j_1)\in T_{k_n,\ell_n}}\sum_{(i_2,j_2)\in T_{k_n,\ell_n}} \!\!\!\!
G_{n,i_1,j_1,i_2,j_2}\\
&\phantom{=}\leq\sum_{(i_1,j_1)\in  T_{k_n,\ell_n}}
\sum_{(i_2,j_2)\in T_{k_n,\ell_n}}\!\!\!\! \Big(
  2M_4L^{(1)}_{n,i_1,j_1,i_2,j_2}\!+\!\! 
 (M_4\!-\!3)^+L^{(2)}_{n,i_1,j_1,i_2,j_2}\Big) \\
&\phantom{==}+
(M_4\!-\!3)^+\!\left(\sum_{i=-\ell_n+1}^{k_n}\sum_{j_1=-i+1}^{\ell_n} 
\sum_{j_2=-i+1}^{\ell_n}\!\!\!\! L^{(3)}_{j_1,j_2}(\alpha_{r_n})+\!\!\!\!\!\!
\sum_{j=-k_n+1}^{\ell_n}\sum_{i_1=-j+1}^{k_n}
\sum_{i_2=-j+1}^{k_n}\!\!\!\! L^{(3)}_{i_1,i_2}(\beta_{r_n})\right),
\end{align*}
where 
\begin{align*}
L^{(1)}_{n,i_1,j_1,i_2,j_2}&:=
\cov\big (X_{i_1-1,j_1}^{(r_n)}-X_{i_1,j_1-1}^{(r_n)},
X_{i_2-1,j_2}^{(r_n)}-X_{i_2,j_2-1}^{(r_n)}\big)^2, \\
L^{(2)}_{n,i_1,j_1,i_2,j_2}&:=\!\!\!\sum_{(u,v)\in 
U_{i_1\land i_2-1,j_1\land
  j_2-1}}(\alpha_{r_n}+\beta_{r_n})^{2(i_1+j_2+i_2+j_2-2-2u-2v)} \\
&\phantom{\qquad \quad \quad =}\times \Big (
{\mathsf P}\big (S^{(\nu_{r_n} )}_{i_1+j_1-1-u-v}=i_1-u\big )
-{\mathsf P}\big (S^{(\nu_{r_n} )}_{i_1+j_1-1-u-v}=i_1-1-u
\big )\Big )^2 \\
&\phantom{\qquad \quad \quad =}\times \Big (
{\mathsf P}\big (S^{(\nu_{r_n} )}_{i_2+j_2-1-u-v}=i_2-u\big )
-{\mathsf P}\big (S^{(\nu_{r_n} )}_{i_2+j_2-1-u-v}=i_2-1-u
\big )\Big )^2, \\
 L^{(3)}_{i_1,i_2}(\nu)&:=\sum_{u=-\infty}^{i_2\land i_2-1}\nu
 ^{2(i_1+i_2-2-2u)}\leq \frac 1{1-\nu^2}, \qquad \qquad 0<|\nu|<1.
\end{align*}
Obviously,
\begin{align*}
&\frac 1{n^4}\sum_{(i_1,j_1)\in T_{k_n,\ell_n}}
\sum_{(i_2,j_2)\in T_{k_n,\ell_n}} L^{(1)}_{n,i_1,j_1,i_2,j_2} \\
& \ =\!\!
\iint\limits _T\!\iint\limits _T\!\!\Big (r_n^{1/2}\cov
\big(Z_{0,1}^{(n)}(s_1,t_1)\!-\! 
Z_{1,0}^{(n)}(s_1,t_1),Z_{0,1}^{(n)}(s_2,t_2)\!-\!Z_{1,0}^{(n)}(s_2,t_2)\big)
\Big )^2{\mathrm d}s_1{\mathrm d}t_1{\mathrm d}s_2{\mathrm d}t_2,
\end{align*}
where due to \eqref{eq:eq1.7}, Propositions \ref{diffbound}, \ref{covlim} 
and Fatou's lemma the right hand side converges to \ $0$ \
as \ $n\to\infty$. \ Furthermore, using Remark \ref{rem:rem4} one can
find an upper bound for \ $L^{(2)}_{n,i_1,j_1,i_2,j_2}$, \ namely
\begin{equation*}
L^{(2)}_{n,i_1,j_1,i_2,j_2}\leq \frac
{D\zeta(5/4)}{(\alpha_{r_n}\beta_{r_n})^3(i_1\lor i_2+j_1\lor 
  j_2)^{1/4}}
\end{equation*}
with some positive constant \ $D$. \  Hence,  
\begin{equation*}
\frac 1{n^4}\sum_{(i_1,j_1)\in T_{k_n,\ell_n}}
L^{(2)}_{n,i_1,j_1,i_2,j_2}
\leq  \frac {20D\zeta(5/4)}{(\alpha_{r_n}\beta_{r_n})^3n^{1/4}}\to 0
\end{equation*}
as \ $n\to\infty$. \ Finally, if  \ $\nu_{r_n}$ \ denotes one of the
sequences \ $\alpha_{r_n}$ \ or \ $\beta_{r_n}$ \ we have 
\begin{equation*}
\lim _{n\to\infty}\frac 1{n^4}\sum_{i=-\ell_n+1}^{k_n}\sum_{j_1=-i+1}^{\ell_n} 
\sum_{j_2=-i+1}^{\ell_n} L^{(3)}_{j_1,j_2}(\nu_{r_n}) \leq 
\lim _{n\to\infty} \frac 1{n(1-\nu_{r_n}^2)}=0,
\end{equation*}
that completes the proof. \proofend

\noindent
\textbf{Proof of Proposition \ref{pro:pro8}.} \ Using the same
techniques as in the proof of
Proposition \ref{pro:pro6} with the help of representation
\eqref{eq:eq6.5} one can show that
\begin{equation*}
\frac 1{n^4}\sum_{m=1}^n
   {\mathsf E}\big((Q_{n,m}-Q_{n,m-1})^4\,\big|\,
{\mathcal F}^n_{m-1}\big)  \stoch  0 \qquad\text{as \ $n\to\infty$.}
\end{equation*}
\space \hfill 
 \proofend

\section*{Acknowledgments}

This research has been supported by  the Hungarian  Scientific
 Research Fund under Grants No. ~OTKA-F046061/2004 and
 ~OTKA-T048544/2005.

\end{document}